\documentstyle[12pt]{article}
\newcommand{\ignore}[1]{}

\begin{document}

\begin{center}{\Large\bf Holey Schr\"oder Designs of Type $\bf 3^n u^1$}

\vskip 16pt

{\sc Dianhua Wu} \\
 School of Mathematics and Statistics\\
 Guangxi Normal University \\
 Guilin 541004, P.R.\ China \\
\texttt{dhwu@gxnu.edu.cn}
\vskip 10pt
{\sc Hantao Zhang} \\
Computer Science Department\\
The University of Iowa\\
Iowa City, IA 52242 U.S.A.\\
\texttt{hantao-zhang@uiowa.edu}

\vskip 0.5 in
\end{center}

\newtheorem{Theorem}{Theorem}[section]
\newtheorem{Definition}[Theorem]{Definition}
\newtheorem{Example}[Theorem]{Example}
\newtheorem{Lemma}[Theorem]{Lemma}
\newtheorem{Construction}[Theorem]{Construction}
\newcommand{\p}{{\noindent \bf Proof}\ \ }
\newcommand{\e}{\hfill $\Box $}

\begin{abstract}
  A holey Schr\"oder design of type
  $h_1^{n_1}h_4^{n_2}\cdots h^{n_k}_k$
  (HSD$(h_1^{n_1}h_4^{n_2}\cdots h^{n_k}_k))$ is equivalent to a frame
  idempotent Schr\"oder quasigroup
  (FISQ$(h_1^{n_1}h_4^{n_2}\cdots h^{n_k}_k))$ of order $n$ with $n_i$
  missing subquasigroups (holes) of order $h_i, 1 \le i \le k$, which
  are disjoint and spanning (i.e., $\sum_{1\le i \le k}n_ih_i =
  n$). The existence of HSD$(h^nu^1)$ for $h=1, 2, 4$ has been known.
  In this paper, we consider the existence of HSD$(3^nu^1)$ and show
  that for $0\le u \le 15$, an HSD$(3^nu^1)$ exists if and only if
  $n(n + 2u -1) \equiv 0~(mod~4)$, $n\ge 4$ and $n\ge 1+2u/3$.  For
  $0 \le u \le n$, an HSD$(3^nu^1)$ exists if and only if
  $n(n + 2u -1) \equiv 0~(mod~4)$ and $n \ge 4$, with possible
  exceptions of $n = 29, 43$. We have also found six new HSDs of type
  $(4^nu^1)$.
\end{abstract}

\section{Introduction}

A {\it quasigroup} is an ordered pair $(Q, * )$, where $Q$ is a set
and ($*$) is a binary operation on $Q$ such that the multiplication
table for $*$ is a Latin square over $Q$. A quasigroup is called {\em
  idempotent} if $x* x = x$ for all $x \in Q$.  If the identity
\begin{equation}
(x* y)* (y* x) = x
\end{equation}
holds for all
$x, y \in Q$, then it is called a {\em Schr\"oder quasigroup}. The
{\em order} of the quasigroup $(Q, *)$ is $|Q|$.

Two quasigroups of the same order are {\it orthogonal} if when the two
corresponding Latin squares are superposed, each symbol in the first
square meets each symbol in the second square exactly once. A
quasigroup (Latin square) is called {\it self-orthogonal} if it is
orthogonal to its transpose. A pair of orthogonal Latin squares
(quasigroups), say $(Q, *)$ and $(Q, \cdot)$, are said to have the
{\em Weisner property} if $x* y = z$ and $x \cdot y = w$ whenever
$z * w = x$ and $z \cdot w = y$ for all $x,y,z,w \in Q$.

Idempotent Schr\"oder quasigroups are associated with other
combinatorial configurations such as a class of edge-colored block
designs with block size 4, triple tournaments and self-orthogonal
Latin squares with the Weisner property (see \cite{cs}, \cite{b},
\cite{lms}, \cite{m} and \cite{weis}).
An idempotent Schr\"oder design of order $v$ over the set $Q$ (i.e., $|Q|=v$)
is equivalent to an
edge-colored design CBD[$G_6;v$] \cite{cs},
which is a partition of the colored edges of a triplicate complete graph
3$K_v$, each edge of $K_v$ receiving one of three
colors, into blocks $[a, b, c, d]$, each containing edges
$\{a, b\}, \{c, d\}$ colored with color 1, edges $\{a, c\}, \{b, d\}$
with color 2, and edges $\{a, d\}, \{b, c\}$ with color 3. If we
define a binary operation ($\cdot $) as
$a\cdot b = c, b\cdot a = d, c \cdot d = a$ and $d\cdot c = b$ from
the block $[a, b, c, d]$ and define $x\cdot x = x$ for every
$x \in Q$, then $(Q, \cdot)$ is an idempotent Schr\"oder quasigroup.  On
the other hand, suppose $(Q, \cdot)$ is an idempotent Schr\"oder
quasigroup.  If $a\cdot b = c$, $b\cdot a = d$, then we must have
$c\cdot d = (a\cdot b)\cdot (b\cdot a) = a$ and
$d\cdot c = (b\cdot a)\cdot (a\cdot b) = b$. So the block
$[a, b, c, d]$ is determined and a CBD[$G_6;v$] can be obtained in
this way.

Let $Q$ be a set and ${\cal H} = \{S_1,S_2,\cdots ,S_k $\} be a set of
subsets  of $Q$ (each subset is called a {\em hole}). A {\em holey
  Schr\"oder quasigroup} having hole set ${\cal H}$ is a triple
$(Q, {\cal H},*)$, which satisfies the following properties:

\begin{enumerate}
\item $(* )$ is a binary operation defined on $Q$, however, when both
  points $a$ and $b$ belong to the same hole, there is no definition
  for $a* b$, and, for any $c\in Q$, $a*c \neq b$ and $c*a \neq b$;
\item no repeated points in each of row or column of the
  multiplication table of $*$;
\item the identity (1) holds when $x$ and $y$ are not in the same hole
  $S_i, 1 \le i \le k$.
\end{enumerate}

If ${\cal H} = \{S_1,S_2,\cdots ,S_k$\} is a partition of $Q$, then a
holey Schr\"oder quasigroup is called {\em frame Schr\"oder
  quasigroup} (FSQ). The {\em type} of the FSQ is defined to be the
multiset \{$|S_i|: 1\le i \le k$\}. We shall use an ``exponential''
notation $s_1^{n_1}s_2^{n_2} \cdots s_t^{n_t}$ to describe the type of
$n_i$ occurrences of $s_i$, $1\le i \le t$, in the multiset. We
briefly denote an FSQ of type $s_1^{n_1} s_2^{n_2}\cdots s_t^{n_t}$ by
FSQ($s_1^{n_1}s_2^{n_2}\cdots s_t^{n_t}$).

Now from an FSQ($s_1^{n_1}s_2^{n_2}\cdots s_t^{n_t}$), we can use
the same method to obtain an edge-colored design which is called a
holey Schr\"oder design and denoted by \linebreak
HSD$(s_1^{n_1}s_2^{n_2}\cdots s_t^{n_t})$.  A {\em  holey Schr\"oder
design} (HSD) is a triple $(X,{\cal H, B})$ which satisfies the following
properties:

\begin{enumerate}
\item ${\cal H}$ is a partition of $X$ into subsets called {\it holes},
\item ${\cal B}$ is a family of 4-subsets of $X$ (called $blocks$) such that a hole
and a block contain at most one common point,
\item the pairs of points in a block $[a, b, c, d]$ are colored as
$\{a, b\}$ and $\{c, d\}$ with color 1, $\{a, c\}$ and $\{b, d\}$
with color 2, and $\{a, d\}$ and $\{b, c\}$ with color 3,
\item every pair of points from distinct holes occurs in 3 blocks with
different colors.
\end{enumerate}

An HSD can be viewed as a generalization of CBD[$G_6;v$]
and is equivalent to a frame self-orthogonal
Latin square (FSOLS) with the Weisner property \cite{bch,lms}.

The {\it type} of the HSD is the multiset $\{|S_i|: S_i \in {\cal H}\}$
and is described by an exponential notation, i.e., we use
$s_1^{n_1}s_2^{n_2} \cdots s_t^{n_t}$ to describe the type of $n_i$
occurrences of $s_i$, $1\le i \le t$, in the multiset.
We denote such an HSD by HSD$(s_1^{n_1}s_2^{n_2} \cdots s_t^{n_t})$.

In particular, an HSD$(1^v)$ is equivalent to an idempotent Schr\"oder
design, and an HSD$(1^{v-n}n^1)$ is also called {\em an incomplete
  Schr\"oder design} (ISD), which has a single hole of size $n$ and
will be denoted by ISD$(v, n)$.

We state the following theorem,
which provides some of the known results on HSDs and ISDs.

\begin{Theorem}\label{known0}{\rm (\cite{lms}, \cite{bwz}, \cite{cs})}

$(a)$ A Schr\"oder quasigroup of order $v$ exists if and only if $v$
$\equiv 0, 1~(mod~4)$ and $v \neq 5$.

$(b)$ An idempotent Schr\"oder quasigroup of order $v$ exists if and
only if $v \equiv 0, 1~(mod~4)$ and $v \neq 5, 9$.  Equivalently, an
ISD$(v, 1)$ exists if and only if $v \equiv 0, 1~(mod~4)$ and
$v \neq 5, 9$.

$(c)$ An ISD$(v, 2)$ exists if and only if $v \ge 7$
and $v \equiv 2, 3~(mod~4)$ except for $v = 10$.

$(d)$ {\rm (\cite{bw, zz})} An HSD$(h^n)$ exists if and only if
$h^2n(n-1)\equiv 0 \ (mod\ 4)$ with exceptions of $(h,n)\in \{(1,5),
(1,9), (2,4)\}$.

$(e)$ {\rm (\cite{HSD2nu1, zz})} For $1 \le u \le 16$, an HSD$(2^nu^1)$
  exists if and only if $n \ge u + 1$, with exception of
  $(n,u) \in \{ (2,1), (3,1), (3,2)\}$, and with possible exceptions
  of $(n, u) \in \{ (7, 5), (7, 6), (11, 9), (11, 10) \}$.
  For $17 \le u \le 4(n-14)/5$, there exists an HSD$(2^nu^1)$.

$(f)$ {\rm (\cite{bz2})} For $0 \le u \le 36$, an HSD$(4^nu^1)$ exists if
and only if $n \ge 4$ and
$0\le u \le 2n -2$,  with possible
exceptions of
$n = 19$ and $u \in \{ 29, 30, 31, 33, 34, 35\}$, or
$n = 22$ and $u \in \{ 33, 34, 35\}$.
For $37 \le u \le 3(n-7)/2$, there exists an HSD$(4^nu^1)$.
\end{Theorem}

The following lemma provides a necessary condition for the existence
of an HSD$(3^n u^1)$.

\begin{Lemma}\label{necessary}
A necessary condition for the existence of an HSD$(3^n u^1)$ is
$n(n + 2u -1) \equiv 0~(mod~4)$, $n\ge 1+2u/3$, and $n\ge 4$.
\end{Lemma}
\p Suppose $(X,{\cal H, B})$ is an HSD$(3^n u^1)$, where
${\cal H} = \{A_1, A_2, ..., A_n, Y\}$, $|A_i| = 3$ for $1 \le i \le n$ and
$|Y| = u$. Let $A = \bigcup_{i=1}^n A_i$, where $|A| = 3n$.
We know that no
block of ${\cal B}$ can contain two or more points from $Y$. Let us
consider pairs of points which have the color 1 in blocks of
${\cal B}$. For any given pair $\{x, y\}$ of points not in the same hole,
there exists one block of the form $[x, y, *, *]$ or
$[*, *, x, y]$ in $\cal B$.
The number of pairs $\{x, y\}$, where $x\in A$ and $y\in Y$, is $3nu$
and each of these pairs must appear exactly in one block. These $3nu$ blocks
also contain $3nu$ pairs $\{a, b\}$ of color 1, where $a, b\in A$ and
$a$ and $b$ cannot be in the same hole $A_i$ for any $i$. Thus,
the total number of such $\{a, b\}$ is $3n(3n-1)/2 - 3n = 3n(3n-3)/2$, and
$3n(3n-3)/2 \ge 3nu$, which implies that $n\ge 1+2u/3$.
Moreover, each block without elements of $Y$ contains two pairs of color 1 from $A$.
Thus, $3n(3n-3)/2 - 3nu$ must be even.
A simple calculation shows that
\[
  (3n(3n-3)/2 - 3nu)/2 = 2n(n - u - 2) + n(n+2u-1)/4
\]
Consequently,
$n(n+2u-1)/4$ must be an integer, that is, $n(n+2u-1)\equiv 0~(mod~4)$.
Computer search shows that no HSD$(3^n u^1)$ exists when $n < 4$.
\e

We say $u$ is {\em feasible} for HSD$(3^n u^1)$ if $u$ satisfies the necessary condition
of Lemma~\ref{necessary}.
According to Lemma~\ref{necessary}, no HSD$(3^n u^1)$ exists if $n \equiv 2~(mod~4)$;
if HSD$(3^n u^1)$ exists, then $0\le u \le 3(n+1)/2$, either $n \equiv 0~(mod~4)$, or
\begin{itemize}
  \item $n \equiv 1~(mod~4)$ and $u$ is even, or
  \item $n \equiv 3~(mod~4)$ and $u$ is odd.
\end{itemize}
The above restrictions on $u$ when $n \not\equiv 0~(mod~4)$ make the
construction of HSD$(3^n u^1)$ more challenge than that of HSD$(2^n u^1)$
or HSD$(4^n u^1)$.

In this paper, we will establish the following results:

\begin{Theorem}\label{main}
  $(a)$ For $0 \le u \le 15$, an HSD$(3^n u^1)$ exists if and only if
  $n(n + 2u -1) \equiv 0~(mod~4)$, $n\ge 1+2u/3$, and
  $n \ge 4$.

  $(b)$ For $0 \le u \le n$,
  an HSD$(3^n u^1)$ exists if and only if
  $n(n + 2u -1) \equiv 0~(mod~4)$ and $n \ge 4$, with possible
  exceptions of $n = 29, 43$.
\end{Theorem}

\section{Construction Methods}

Among other things, we shall make use of some known constructions
relating HSDs and other designs.
To construct HSDs directly, sometimes we can use {\it starter
blocks}. Suppose the block set ${\cal B}$ of an HSD is closed under
the action of some Abelian group $G$, then we are able to list only
part of the blocks (starter or base blocks) which determines the
structure of the HSD. To check the starter blocks, we need only
calculate whether the differences $\pm (x - y)$ from all pairs $\{x,
y\}$ with color $i$ in the starter blocks are precisely $G\setminus
S$ for $1\le i \le 3$, where $S$ is the set of the differences of
the holes.  We can also attach some infinite points to an Abelian
group $G$.  When the group acts on the blocks, the infinite points
remain fixed. In the following example $x_i$ is an infinite
point.

\begin{Example}\label{ex3711}\rm An HSD$(3^7 1^1)$.

Points $X$: $Z_{21}  \cup \{x\}$

Holes $\cal H$:  $\{\{ i, i+7, i+14 \} : 0 \le i \le 7 \} \cup \{x\}$

Starter blocks:
$[0,  1,  5, x], [0,  2, 12,  1], [0,  3, 18,  9], [0,  4,  2,  8], [0,  5, 10, 18]$.

In this example, the entire set of blocks are developed from the
starter blocks by repeatedly adding $1~(mod~21)$ to each point of the
starter blocks; the infinite points are unchanged for addition.
The result is a set of $105$ blocks.
\end{Example}

The above idea of starter blocks can be also generalized: Instead of
adding 1 to each point of the starter blocks, we may add 2 or more to
develop the block set; we refer this as the +2 method.

\begin{Example}\label{ex3821}\rm An HSD$(3^8 2^1)$.

Points $X$: $Z_{24}  \cup \{x_1, x_2\}$

Holes $\cal H$:  $\{\{ i, i+8, i+16 \} : 0 \le i \le 7 \} \cup \{x_1, x_2 \}$

Starter blocks:
$\bf [0,  9, 12, 21]$,
$\bf [0, 11, 23, 12]$,
$\bf [0, 12,  3, 15]$,
$[0,  1, 20,  3], [0,  2, 19, 13]$, $[0,  3, 18, 14], [0,  5,  9, 23]$,
$[0, 10,  5,  7]$,
$[0, 15, 14, 11], [0, 17, 15, x_1], [0, 18, 22, x_2]$, $[0, 19,  1,  2]$,
$[1,  5, 15, x_2]$, $[1, 12, 18, x_1]$.

In this example, the entire set of blocks are developed from the
starters by repeatedly adding $2~(mod~24)$ to each point of the
starters; the infinite points are unchanged for addition.  The
result is a set of $14\times 12 = 168$ blocks. However,
each of the first three starters will generate only $6$ instead of $12$ distinct blocks by
the +2 method. For instance, $[0, 12, 3, 15]$ generates the following six
distinct blocks:
\[ [0, 12, 3, 15], [2, 14, 5, 17], [4, 16, 7, 19], [6, 18, 9, 21], [8, 20, 11, 23], [10, 22, 13, 1],
\]
because $[12, 0, 15, 3]$ is equivalent to $[0, 12, 3, 15]$, etc.
Similarly, $[0,  9, 12, 21]$ and $[0, 11, 23, 12]$ can each generate $6$ distinct blocks.
In other words,
these three starters are invariant under the operation $+12~(mod~ 24)$ and
use so-called {\em short orbit} to generate blocks.
The total number of distinct blocks for this example is $150$ ($3\times 6+11\times 12$).
\end{Example}

The following materials come from {\rm \cite{bz,bz2}}.
These recursive constructions of HSDs are commonly used
in other block designs.

Another class of designs related to HSDs is {\it group divisible
design} (GDD) \cite{bsh}. A GDD is a 4-tuple $(X, {\cal G, B, \lambda}$) which
satisfies the following properties:
\begin{enumerate}
\item ${\cal G}$ is a partition of $X$ into subsets called {\it groups},
\item ${\cal B}$ is a family of subsets of $X$ (called {\it blocks}) such that
a group and a block contain at most one common point,
\item every pair of points from distinct groups occurs in exactly $\lambda$
blocks.
\end{enumerate}
The {\it type} of the GDD is the multiset $\{ |G|: G \in {\cal
  G}\}$. We also use the notation GD($K,M;\lambda $) to denote the GDD
when its block sizes belong to $K$ and group sizes belong to $M$.  If
$M =\{1\}$, then the GDD becomes a {\em pairwise balanced design} (or
{\em PBD}) \cite{w,mg}. If $K = \{k\}$, $M = \{n\}$ and with the type $n^k$, then
the GDD becomes TD$(k,n)$.  It is well known that the existence of a
TD$(k,n)$ is equivalent to the existence of $k - 2$ MOLS$(n)$ \cite{abcd}. It is
easy to see that if we erase the colors in the blocks, the HSD becomes
a GDD with block size 4 and $\lambda = 3$. 

\begin{Theorem}\label{TD}
$(a)$
There exists a TD$(6, m)$
for $m \ge 5$ and $m \not \in \{6, 10, 14, 18, 22\}$,

$(b)$ There exists a TD$(r, p)$ if $p$ is a prime power and $r\le p+1$;
\end{Theorem}

The following construction comes from the weighting construction of
GDDs \cite{w}.

\begin{Construction}\label{Weight} {\rm \cite{w}}
  $($Weighting$)$ Suppose $(X,{\cal H, B})$ is a GDD with
  $\lambda = 1$ and let $w: X \mapsto Z^+\cup\{0\}$. Suppose there
  exist HSDs of type $\{w(x): x \in B\}$ for every $B \in {\cal
    B}$. Then there exists an HSD of type
  $\{\sum_{x\in H} w(x): H \in {\cal H} \}$.
\end{Construction}

The following lemma uses the weighting construction.

\begin{Lemma}\label{GDD3nu1} {\rm (\cite{bz2})}
  There exists an HSD$(3^nu^1)$ if one of the following is true:

~~~~  $(a)$ $n \equiv 0 ~(mod~ 4)$ and $u \equiv 0~(mod~3)$, $0 \le u \le (3n - 6)/2$;

~~~~  $(b)$ $n \equiv 1~ (mod~ 4)$ and $u \equiv 0~(mod~6)$, $0 \le u \le (3n - 3)/2$; or

~~~~  $(c)$ $n \equiv 3~ (mod~ 4)$ and $u \equiv 3~(mod~6)$, $0 < u \le (3n - 3)/2$.
\end{Lemma}
\p There exist 4-GDDs of the same type (see, for example,
\cite{g}). So we can give all points of this GDD weight one to get
the desired HSD$(3^{n} u^1)$.  \e

We need HSD$(3^n u^1)$ of small $n$ for
the weighting construction using TDs.

\begin{Lemma}\label{3nsmall}
  There exists an HSD$(3^n u^1)$ if
  $n=4$ and $0\le u \le 4$, or $n=5$ and $u\in \{0, 2, 4, 6\}$.
\end{Lemma}
\p Theorem~\ref{known0}$(d)$ covers the case when $u=0$.
$(n,u)=(4,3)$ is equivalent to $(n,u)=(5,0)$.
$(n,u)=(5,6)$ is covered by Construction\ref{GDD3nu1}$(b)$.
The other HSDs can be found in \cite{bz2}.
\e

\begin{Lemma}\label{TD6mtk1u1}
  If there exists a TD$(6, m)$, then there exists an
  HSD$((3m)^4(3k)^1 u^1)$, where $0 \le k \le m$, $u$ is even and
  $0 \le u\le 4m$.
\end{Lemma}

\p  Give weight 3 to each point of the first four groups of a TD$(6,
m)$. Give a weight of $3$ to $k$ points of the fifth group and
weight 0 to the remaining points of this group. Give a weight of 0,
2, 4 to each point of the sixth group so that the
total weight is $u$. By using HSDs of types
$3^n k^1$ for $n = 4, 5$ and $k = 0, 2, 4$ from
Lemma~\ref{3nsmall}, we obtain the desired HSD by
Construction~\ref{Weight}. \e

\begin{Construction}\label{Multiply}{\rm (\cite{bz2})}
  Suppose there exists an HSD$(h_1^{n_1} h_2^{n_2} \cdots h_k^{n_k})$,
  then there exists an
  HSD$((mh_1)^{n_1}(mh_2)^{n_2}\cdots (mh_k)^{n_k})$, where
  $m \neq 2,6$.
\end{Construction}

\begin{Lemma}\label{9n4u1}
  An HSD$(9^{4} u^1)$ exists for $0\le u \le 15$.
\end{Lemma}
\p If $u \equiv 0~(mod~3)$, let $u=3t$, where $0\le t \le 4$, and
 there exists an HSD$(3^4 t^1)$ by Lemma~\ref{3nsmall}.
 To this HSD we apply Construction~\ref{Multiply} with $m=3$
 and obtain the desired HSD.
For $u \not\equiv 0~(mod~3)$ and $1 \le u \le 13$,
an HSD$(9^4 u^1)$ is provided in Appendix C1. \e

\begin{Lemma}\label{9n5u1}
  An HSD$(9^{5} u^1)$ exists for even $u$ and $0\le u \le 18$.
\end{Lemma}
\p If $u \equiv 0~(mod~3)$, let $u=3t$, where $t$ is even and $0\le t \le 6$, and
 there exists an HSD$(3^5 t^1)$ by Lemma~\ref{3nsmall}.
 To this HSD we apply Construction~\ref{Multiply} with $m=3$
 and obtain the desired HSD.
For $u \not\equiv 0~(mod~3)$ and $2 \le u \le 16$,
an HSD$(9^5 u^1)$ is provided in Appendix C2.
\e

\begin{Lemma}\label{9n9u1}
  An HSD$(9^{9} u^1)$ exists for even $u$ and $8\le u \le 36$.
\end{Lemma}
\p For $u=8, 10, 14, 16$, the HSDs can be found in Appendix C3.
For $u=12$, take HSD$(3^9 4^1)$ from Appendix A3, apply Construction~\ref{Multiply} with
$m=3$. For $18\le u \le 36$,
we start with a TD$(10, 9)$
and give all the
points in all the groups except the last group a weight of one. In the
last group we can give $k$ points a weight of 4 and the rest a weight
of 2 so that the total weight $u$ of the group satisfies $u = 4k + 2(9-k)$
and $18 \le u \le 36$.  Here we need HSDs of types $1^9 u^1$ for $u=2,4$
which exist \cite{bz2}. The resulting design is the desired result.  \e

With the existence result on HSDs, proving
the existence of HSD$(h^n u^1)$ becomes much easier due to the next
construction, which may be called ``filling in holes''. It is used
commonly in constructing designs.

\begin{Construction}\label{fillinghole} {\rm (\cite{bz2})}
$(a)$  If there exist an HSD$(h^sv^1)$
and an HSD$((hs)^mw^1)$, then there exists an HSD$(h^{sm}(w+v)^1)$.

$(b)$ If there exist an HSD$(h^sv^1)$, an HSD$(h^tv^1)$,
and an HSD$((hs)^m (ht)^1 w^1)$, then there exists an HSD$(h^{sm+t}(w+v)^1)$.
\end{Construction}
\p  In both cases, we add $v$ points to HSD$((hs)^m w^1)$ or
 HSD$((hs)^m (ht)^1 w^1)$,
 fill in the $m$ holes of size $hs$ with an HSD$(h^s v^1)$,
 and fill in the one hole of size $ht$ with an HSD$(h^t v^1)$ in case $(b)$.
\e

For $n=15$, $u \not\equiv 3~(mod~6)$ and $5\le u \le 19$, by Lemma~\ref{9n5u1},
there exists an HSD$(9^5 (u-3)^1)$. To this HSD,
we apply Construction~\ref{fillinghole}$(a)$ with $h=s=v=3$, $m=5$, and $w=u-3$.

\begin{Construction}\label{great} {\rm (\cite{bz2})}
Suppose there exists an FSOLS$(h_1^{n_1}h_2^{n_2} \cdots h_k^{n_k})$,
then there exists an HSD$((4h_1)^{n_1}(4h_2)^{n_2} \cdots (4h_k)^{n_k})$.
\end{Construction}

To use the above construction, we need the following result.

\begin{Theorem}\label{fsols} {\rm (\cite{xz,xzz})}
  An FSOLS$(3^m k^1)$ exists if and only if $m \ge 1+2k/3$, $m \ge 4$,
  $k \ge 0$, except possibly
\footnote{FSOLS of types $(3^6 7^1)$, $(3^{14} 19^1)$, and $(3^{18} 25^1)$
were found by the second author.}
  $(m, k) = (22, 31)$.
\end{Theorem}

\begin{Lemma}\label{12m4k1} {\rm (\cite{bz2})}
  An HSD$(12^m (4k)^1)$ exists for all $m \ge 1+2k/3$, $m \ge 4$, and
  $k \ge 0$, except possibly
  $(m, k) = (22, 31)$.
\end{Lemma}
\p By Theorem~\ref{fsols}, there exists an FSOLS$(3^m k^1)$. Applying
Construction~\ref{great} to this FSOLS, we obtain an HSD$(12^m (4k)^1)$ with
the possible exception listed in the lemma.  \e

\section{Preliminary Results of HSD$(3^n u^1)$}

Using the constructions presented above, we can establish the following results.

\begin{Lemma}\label{3n12u1}
For $n=8, 12$, an HSD$(3^{n} u^1)$ exists for all feasible value of $u$.
\end{Lemma}
\p Lemma~\ref{necessary} tells that $0 \le u \le 10$ for $n=8$ and
$0 \le u \le 16$ for $n=12$.
Lemma~\ref{GDD3nu1}$(a)$ covers the cases when $u \equiv 0~(mod~3)$.
For $n=8$ and $u=2$, see Example~\ref{ex3821}; the other cases of
$u=1, 2$ are provided in Appendix~A1 and A2, respectively.

For $n=8$ and $u=4,5,7,8, 10$, the designs are provided in Appendix
A3, A4, A5, A6, and A7, respectively.

For $n=12$, $u \not\equiv 0~(mod~3)$ and $4\le u \le 16$, by Lemma~\ref{9n4u1},
there exists an HSD$(9^4 (u-3)^1)$. To this HSD,
we apply Construction~\ref{fillinghole}$(a)$ with $h=s=v=3$, $m=4$, and $w=u-3$.
\e

\begin{Lemma}\label{lem3n11}
An HSD$(3^n 1^1)$ exists for all $n \ge 4$ and $n \equiv 0, 3~(mod~4)$.
\end{Lemma}
\p For $n \equiv 0~(mod~4)$,
the cases of $n = 4, 8, 12$ are provided in Lemma~\ref{3nsmall},
Appendix A1, and Lemma~\ref{3n12u1}, respectively.
For $n \ge 16$, let $n = 4m$ for some $m \ge 4$.  By
Theorem~\ref{known0}$(d)$, there exists an HSD$(12^m)$,
and we apply Construction~\ref{fillinghole}$(a)$ with $h=3$, $s=4$, $v=1$, and $w=0$.

For $n \equiv 3~(mod~4)$,
see Example~\ref{ex3711} for $n=7$;
the cases of $n = 11, 15$ are provided in Appendix A1.
For $n = 19$,
we apply Construction~\ref{fillinghole}$(a)$ with $h=s=v=3$, $m=5$, and $w=10$,
to the HSD$(9^5 10^1)$ in Appendix C2 to obtain an HSD$(3^{15} 13^1)$,
and then fill in the hole of size $13$ with an HSD$(3^4 1^1)$ (Lemma~\ref{3nsmall}).

For $n = 23$, we apply
Construction~\ref{Multiply} with $m = 3$ to the HSD$(4^4 6^1)$ from
Theorem~\ref{known0}$(d)$ to obtain an HSD$(12^4 18^1)$, and then
apply Construction~\ref{fillinghole}$(a)$ with $h=3$, $s=4$, $m=4$, $v=4$, and $w=18$,
to obtain an HSD$(3^{23} 22^1)$.
Finally, we fill in the hole of size $22$ with an HSD$(3^7 1^1)$.

For $n \ge 27$, let $n = 4m+7$ for some
$m \ge 5$.  By Lemma~\ref{12m4k1} with $k =5$, there exists an
HSD$(12^m 20^1)$ and
we apply Construction~\ref{fillinghole}$(a)$ with $h=3$, $s=4$, $v=2$, and $w=20$,
to obtain an HSD$(3^{4m} 22^1)$.
Fill in the hole of size $22$ with an HSD$(3^7 1^1)$, we obtain an
HSD$(3^{4m+7} 1^1)$.  \e

\begin{Lemma}\label{3n15u1}
For $n = 7, 11, 15$, an HSD$(3^{n} u^1)$ exists for all feasible value of $u$.
\end{Lemma}
\p Lemma~\ref{necessary} tells that $u$ is odd, $0 \le u \le 9$ for $n=5$,
$0 \le u \le 15$ for $n=11$, and $0 \le u \le 21$ for $n=15$.
Lemma~\ref{GDD3nu1}$(c)$ covers the cases when $u \equiv 3~(mod~6)$.  Lemma~\ref{lem3n11}
covers the cases when $u=1$.

For $n=7, 11$ and $u = 5, 7$, the designs are provided in Appendix A4 and A5.
For $n=11$ and $u = 11, 13$, the designs are provided in Appendix A8 and A9.

For $n=15$, $u \not\equiv 3~(mod~6)$ and $5\le u \le 19$, by Lemma~\ref{9n5u1},
there exists an HSD$(9^5 (u-3)^1)$. To this HSD,
we apply Construction~\ref{fillinghole}$(a)$ with $h=s=v=3$, $m=5$, and $w=u-3$.
\e

\begin{Lemma}\label{3n27u1}
An HSD$(3^{27} u^1)$ exists for odd $u$ and $11\le u \le 39$.
\end{Lemma}
\p
By Lemma~\ref{9n9u1},
there exists an HSD$(9^9 (u-3)^1)$. To this HSD,
we apply Construction~\ref{fillinghole}$(a)$ with $h=s=v=3$, $m=9$, and $w=u-3$.
\e

\begin{Theorem}\label{TD6m4m}
  An HSD$(3^{4m}u^1)$ exists if and only if $m \ge 1$ and $0 \le u\le 6m-2$.
\end{Theorem}
\p The necessary condition comes from Lemma~\ref{necessary}.
For $m=1,2,3$, see Lemmas~\ref{3nsmall} and \ref{3n12u1}.
For $m\ge 4$,
there exists an HSD$(12^m (4k)^1)$ for all
$m \ge 4$ and $0 \le k \le (3m-3)/2$
by Lemma~\ref{12m4k1}.
Since there exists
also an HSD$(3^4 t^1)$ for $0\le t \le 4$ (Lemma~\ref{3nsmall}),
we may add $t$ points to HSD$(12^m (4k)^1)$ and
fill the holes of size $12$ with an HSD$(3^4 t^1)$. The result is
an HSD$(3^{4m}u^1)$, where $u=4k+t$, $0\le u \le 6m-2$, and
$(m, k) \neq (22, 31)$.

Now let us consider the cases of $u=4k+t$ missed by the possible
exception of Lemma~\ref{12m4k1}, i,e. $(m, k) = (22, 31)$. Assuming $k=31$.
For $t=0$, the value of $u=4k+0$
can be obtained from $u=4k'+4$, where $k' = k-1 = 30$.
For $t=4$, the value of $u=4k+4$ can be obtained from $u=4k'+0$, where $k'=k+1=32$.
When $t=2$, $u = 4k+2 = 126 \equiv 0~(mod~3)$,
and in this case there exists an HSD$(3^{4m} u^1)$ by Lemma~\ref{GDD3nu1}$(a)$.
Thus the missed cases of HSD$(3^{n} u^1)$ due to $(m, k) \neq (22, 31)$, $n=4m$, and $u=4k+t$
are: $n = 88$ and $u = 125, 127$.
Now we use
Lemma~\ref{GDD3nu1}$(c)$ to obtain an HSD$(3^{11} 15^1)$.
Apply Construction~\ref{Multiply} (with $m=8$) to this HSD,
we obtain an HSD$(24^{11} 120^1)$.  To this HSD,
we apply Construction~\ref{fillinghole}$(a)$ with $h=3$, $s=8$, $m=7,11$,
$v=t=5,7$, and $w=120$, and obtain an HSD$(3^{88} s^1)$ where
$s \in \{ 125, 127\}$.
The needed HSD$(3^8 t^1)$, $t=5,7$, are available by Lemma~\ref{3n12u1}.
\e

\begin{Lemma}\label{lem3n21}
  An HSD$(3^n 2^1)$ exists for all $n \ge 4$ and
  $n \equiv 0, 1~(mod~4)$.
\end{Lemma}
\p For $n \equiv 0~(mod~4)$, Theorem~\ref{TD6m4m} applies.
For $n \equiv 1~(mod~4)$, see Lemma~\ref{3nsmall} for $n = 5$;
the cases of $n = 9, 13, 17$ are given in Appendix A2.
For $n \ge 21$, let $n = 4m+5$ for some $m \ge 4$.  By Lemma~\ref{12m4k1}
with $k = 4$, there exists an HSD$(12^m 16^1)$.  Adding one point
to the HSD and filling the holes of size 12 with an HSD$(3^4 1^1)$ and
the hole of size $17$ with an HSD$(3^5 2^1)$, we obtain an
HSD$(3^{4m+5} 2^1)$.  \e

\begin{Lemma}\label{lem3n41}
  An HSD$(3^n 4^1)$ exists for all $n \ge 4$ and
  $n \equiv 0, 1~(mod~4)$.
\end{Lemma}
\p For $n \equiv 0~(mod~4)$, Theorem~\ref{TD6m4m} applies.
For $n \equiv 1~(mod~4)$,
see Lemma~\ref{3nsmall} for $n=5$;
the cases of $n = 9, 13, 17$ are provided in Appendix A3.
For $n \ge 21$,
let $n = 4m+5$ for some
$m \ge 4$.  Adding three points to the HSD$(12^m 16^1)$ from
Lemm~\ref{12m4k1} and filling the holes of size 12 with an HSD$(3^5)$
and the hole of size $19$ with an HSD$(3^5 4^1)$, we obtain an
HSD$(3^{4m+5} 4^1)$.  \e

\begin{Lemma}\label{lem3n51}
  An HSD$(3^n 5^1)$ exists for all $n \ge 5$ and $n \equiv 0, 3~(mod~4)$.
\end{Lemma}
\p For $n \equiv 0~(mod~4)$, Theorem~\ref{TD6m4m} applies.
For $n \equiv 3~(mod~4)$,
the cases of $n=7, 11, 15$ are covered by Lemma~\ref{3n15u1};
the cases of $n = 19, 23$ are provided in Appendix A4.
For $n \ge 27$, let $n = 4m+7$ for some
$m \ge 5$.  By Lemma~\ref{12m4k1} with $k = 6$, there exists an
HSD$(12^m24^1)$ for $m \ge 5$. Adding two points to this HSD and
filling the holes of size $12$ with an HSD$(3^42^1)$ and the hole of
size $26$ with an HSD$(3^7 5^1)$, we obtain an HSD$(3^{4m+7} 5^1)$.  \e

\begin{Lemma}\label{lem3n71}
  An HSD$(3^n 7^1)$ exists for all $n \ge 5$ and $n \equiv 0, 3~(mod~4)$.
\end{Lemma}
\p
For $n \equiv 0~(mod~4)$, Theorem~\ref{TD6m4m} applies.
For $n \equiv 3~(mod~4)$, the cases of $n = 7, 11, 19, 23$ are provided in Appendix A6.
For $n=15$, see Lemma~\ref{3n15u1}.
For $n \ge 27$, let $n = 4m+7$ for
$m \ge 5$.  By Lemma~\ref{12m4k1} with $k = 6$, there exists an
HSD$(12^m24^1)$ for $m \ge 5$. Adding 4 points to this HSD and
filling the holes of size 12 with an HSD$(3^44^1)$ and the hole of
size $28$ with an HSD$(3^7 7^1)$, we obtain an HSD$(3^{4m+7} 7^1)$.  \e

\begin{Lemma}\label{lem3n81}
  An HSD$(3^n 8^1)$ exists for all $n \ge 8$ and $n \equiv 0, 1~(mod~ 4)$.
\end{Lemma}
\p
For $n \equiv 0~(mod~4)$, Theorem~\ref{TD6m4m} applies.
For $n \equiv 1~(mod~4)$, the cases of
$n \in \{ 9, 13, 17 \}$ are given in Appendix A6.

For $n = 9, 13, 17, 29$, using the starter blocks from Appendix A6, we
obtain the desired HSDs.  For $n=21, 33$, let $t=4m+1$ for $m=5, 8$.
Since there exists a TD$(6, m)$ for $m=5, 8$, by
Lemma~\ref{TD6mtk1u1}, there exists an HSD$((3m)^4 3^1 8^1)$ with
$m=5, 8$, $k=1$, and $u=8$.  Fill the holes of size $3m$ with
HSD$(3^m)$, we obtain an HSD$(3^{4m+1}8^1)$.

For $n=25$, we use the following starter blocks to generate
an HSD$(5^5 2^1)$ by adding $1~(mod~25)$:
\[[0, 1, 2,x_2], [0, 2, 8,16], [0, 3, 6,12], [0, 7,24,11], [0, 9, 7,21], [0,21, 3, x_1]\]
We then apply
Construction~\ref{Multiply} to this HSD with $m=3$ to obtain an
HSD$(15^5 6^1)$. Adding 2 points to this HSD and filling the holes of
size $15$ with HSDs of type $(3^5 2^1)$, we obtain the desired HSD.

For $n \equiv 1~(mod~ 4)$ and $n \ge 37$, let $n = 4m+9$, where
$m \ge 7$.  By Lemma~\ref{12m4k1} with $k = 8$, there exists an
HSD$(12^m32^1)$ for $m \ge 7$. Adjoining $3$ points to this HSD and
filling the holes of size $12$ with an HSD$(3^5)$ and the hole of
size $35$ with an HSD$(3^9 8^1)$, we obtain an HSD$(3^{4m+9} 8^1)$.
\e

\begin{Lemma}\label{TD6m4m1}
  An HSD$(3^{4m+1}u^1)$ exists if $m \ge 5$,
  $u$ is even and $0 \le u\le 4m+1$, with possible exception of $m=7$.
\end{Lemma}
\p Note that $u \le 4m$ is equivalent to $u\le 4m+1$, since $u$ is even.

Let us first consider $m \equiv 0, 1~(mod~4)$.
Since there exists a TD$(6, m)$,
by Lemma~\ref{TD6mtk1u1} with $k=1$,
there exists an HSD$((3m)^4 3^1 u^1)$, where $u$ is even and $0 \le u\le 4m$.
Filling the holes of size $3m$ with HSD$(3^m)$ (Theorem~\ref{known0}$(d)$),
we obtain an HSD$(3^{4m+1}u^1)$.

Now we consider $m \equiv 2~(mod~4)$.
Applying Lemma~\ref{TD6mtk1u1} with TD$(6, m-1)$ (Theorem~\ref{TD}) and $k=5$,
we obtain an HSD$((3(m-1))^4 15^1 u^1)$, where $u$ is even and $0 \le u\le 4(m-1)$.
Filling the holes of size $3(m-1)$ with HSD$(3^{m-1})$
and the hole of size $15$ with HSD$(3^5)$ (Theorem~\ref{known0}$(d)$),
we obtain an HSD$(3^{4m+1}u^1)$, where $u$ is even and $0 \le u\le 4(m-1)$.
For $4(m-1) \le u \le 4m$, we add $4$ points to
HSD$((3(m-1))^4 15^1 u^1)$ before filling the holes with
HSDs of types $(3^{m-1} 4^1)$ and $(3^5 4^1)$ (Lemma~\ref{lem3n41}).
The result is HSD$(3^{4m+1}u^1)$, where $u$ is even and $4 \le u\le 4m+1$.

Finally, for $m \equiv 3~(mod~4)$,
and $m\ge 11$, applying
Lemma~\ref{TD6mtk1u1} with TD$(6, m-2)$ (Theorem~\ref{TD}) and $k=9$,
there exists an HSD$((3(m-2))^4 27^1 u^1)$, where
$u$ is even and $0 \le u\le 4(m-2)$.
Filling the holes of size $3(m-2)$ with HSD$(3^{m-2})$
and the hole of size $27$ with an HSD$(3^9)$ (Theorem~\ref{known0}$(d)$),
we obtain an HSD$(3^{4m+1}u^1)$, where $0 \le u\le 4(m-2)$.
For $4(m-2) < u \le 4m$, we may add $8$ points to the HSD$((3(m-2))^4 27^1 u^1)$,
and fill the holes with HSDs of types $(3^{m-2} 8^1)$ and $(3^9 8^1)$, which
exist by Lemma~\ref{lem3n81}. This leaves $m=7$ as the only possible exception.
\e

\begin{Lemma}\label{TD6m4m3}
  An HSD$(3^{4m+3}u^1)$ exists if $m \ge 7$,
  $u$ is odd and $0 \le u\le 4m+3$,
  with possible exception of $m = 10$.
\end{Lemma}
\p Note that, since $u$ is odd, $0 \le u$ is equivalent to $1 \le u$.
Because of Lemma~\ref{lem3n11}, we will focus on $3 \le u \le 4m+3$.

Let us first consider $m \equiv 0, 3~(mod~4)$ and $m\ge 7$,
Under the assumption, there exist a TD$(6,m)$ and an HSD$(3^{m+1})$.
By Lemma~\ref{TD6mtk1u1} with $k=3$,
there exists an HSD$((3m)^4 9^1 u^1)$, where $u$ is even and $0 \le u\le 4m$.
We add 3 points to this HSD, then fill all the holes except the hole of size $u$
with HSDs of types $3^{m+1}$ and $3^4$, and obtain an HSD$(3^{4m+3} (u+3)^1)$,
where $3\le u+3 \le 4m+3$.

Now let us consider $m \equiv 1~(mod~4)$, $m\ge 8$.
Since $m-1 \equiv 0~(mod~4)$ and $m-1\ge 7$,
there exists a TD$(6, m-1)$.
By Lemma~\ref{TD6mtk1u1} with $k=7$,
there exists an HSD$((3(m-1))^4 21^1 u^1)$, where $u$ is even and $0 \le u\le 4(m-1)$.
Adding $t=1, 7$ to this HSD, we
fill the holes of size $3(m-1)$ with HSD$(3^{m-1} t^1)$,
and the hole of size $21$ with an HSD$(3^7 t^1)$
and obtain an HSD$(3^{4(m-1)+7} (u+t)^1)$, or HSD$(3^{4m+3} (u+t)^1)$
with odd $u+t$ and $1 \le u+t\le 4m+3$ (Construction~\ref{fillinghole}$(b)$).
The needed HSDs of type $(3^s t^1)$, where $s\in \{ 7, m-1\}$,
are provided in Lemma~\ref{lem3n11} for $t=1$
and Lemma~\ref{lem3n71} for $t=7$.

Finally, for $m \equiv 2~(mod~4)$ and $m\ge 14$,
since $m-2 \equiv 0~(mod~4)$ and $m-2\ge 12$,
there exists a TD$(6, m-2)$.
By Lemma~\ref{TD6mtk1u1} with $k=11$,
there exists an HSD$((3(m-2))^4 33^1 u^1)$, where $u$ is even and $0 \le u\le 4(m-2)$.
Adding $t = 1$ or $11$ to this HSD, we then
fill the holes of size $3(m-2)$ with HSD$(3^{m-2} t^1)$
and the hole of size $33$ with HSD$(3^{11} t^1)$, and
obtain an HSD$(3^{4m+3}(u+t)^1)$, where $u$ is odd and $1 \le u\le 4m+3$
(Construction~\ref{fillinghole}$(b)$).
The needed HSDs of types $(3^s 1^1)$ are provided in Lemma~\ref{lem3n11};
the needed HSDs of type $(3^s 11^1)$ are provided in
Appendix A8 when $s=11$ and by Theorem~\ref{TD6m4m} when $s \equiv 0~(mod~4)$.
This leaves $m=10$ as the only possible exception.
\e

For the possible exception of $m=10$ in the above lemma, we have the following result for
$n=4m+3 = 43$.

\begin{Lemma}\label{3n43}
  An HSD$(3^{43} u^1)$ exists for odd $u$ and $1\le u \le 15$.
\end{Lemma}
\p For $u=1,5, 7$, see Lemmas~\ref{lem3n11}, \ref{lem3n51}, and \ref{lem3n71}.
For $u=3,9, 15$, see Lemma~\ref{GDD3nu1}.
For $u=11, 13$, by Theorem~\ref{TD6m4m} with $m=8$,
we have an HSD$(3^{32}s^1)$, where $0\le s\le 46$.
Using this HSD, let $s=33+t$, where $t=11,13$, and fill the hole of $s$
by an HSD$(3^{11}t^1)$ (available by Lemma~\ref{3n15u1}) to
obtain the desired HSD. \e

\begin{Lemma}\label{lem3n101}
  An HSD$(3^n 10^1)$ exists for all $n \ge 8$ and
  $n \equiv 0, 1~(mod~4)$. 
\end{Lemma}
\p
For $n \equiv 0~(mod~4)$, Theorem~\ref{TD6m4m} applies.
For $n \equiv 1~(mod~4)$,
the cases of $n \in \{ 9, 13, 17, 29 \}$ are provided in Appendix A7.
Lemma~\ref{TD6m4m1} covers the cases when $n \ge 21$ and $n\neq 29$.
\e

\begin{Lemma}\label{lem3n111}
  An HSD$(3^n 11^1)$ exists for all $n \ge 9$ and $n \equiv 0, 3~(mod~4)$.
\end{Lemma}
\p For $n \equiv 0~(mod~4)$,
Theorem~\ref{TD6m4m} applies.
For $n \equiv 3~(mod~4)$,
the cases of $n = 11, 15$ are covered by Lemma~\ref{3n15u1}.
For $n=19, 23$, see Appendix A8.
For $n=27, 43$, Lemmas~\ref{3n27u1} and \ref{3n43} apply.
For $n \ge 31$ and $n \neq 43$, Lemma~\ref{TD6m4m3} applies.
\e

\begin{Lemma}\label{lem3n131}
  An HSD$(3^n 13^1)$ exists for all $n \ge 10$ and $n \equiv 0, 3~(mod~4)$.
\end{Lemma}
\p For $n \equiv 0~(mod~4)$,
Theorem~\ref{TD6m4m} applies.
For $n \equiv 3~(mod~4)$,
the cases of $n = 19, 23$ are provided in Appendix A9.
For $n=11, 15$, see Lemma~\ref{3n15u1}.
For $n=27, 43$, see Lemmas~\ref{3n27u1} and \ref{3n43}.
For $n \ge 31$ and $n \neq 43$, Lemma~\ref{TD6m4m3} applies.
\e

\begin{Lemma}\label{lem3n141}
  An HSD$(3^n 14^1)$ exists for all $n \ge 11$ and
  $n \equiv 0, 1~(mod~4)$.
\end{Lemma}
\p
For $n \equiv 0~(mod~4)$, Theorem~\ref{TD6m4m} applies.
For $n \equiv 1~(mod~4)$,
the cases of $n \in \{ 13, 17, 29 \}$ are provided in Appendix A10.
Lemma~\ref{TD6m4m1} covers the cases when $n \ge 21$ and $n\neq 29$.
\e

\begin{Lemma}\label{3n17u1}
For $n = 13, 17$, an HSD$(3^{n} u^1)$ exists for all feasible values of $u$.
\end{Lemma}
\p Lemma~\ref{necessary} tells that $u$ is even, $0 \le u \le 18$ for $n=13$.
and $0 \le u \le 24$ for $n=17$.
Lemma~\ref{GDD3nu1}$(b)$ covers the cases when $u$ is even and
$u \equiv 0~(mod~3)$.  For $u=2, 4, 8$, $10$, and $14$, the designs are
covered by Lemmas~\ref{lem3n21}, \ref{lem3n41}, \ref{lem3n81},
\ref{lem3n101}, and \ref{lem3n141}, respectively.  For $(n,u) = (13, 16)$, the starter blocks
are $(+1~ mod~ 39)$:
\[\begin{array}{l}
~   [0,  1,  2, x_{16}],
    [0,  2,  4, x_{15}],
    [0,  3,  6, x_{14}],
    [0,  4,  8, x_{13}],
    [0,  5, 10, x_{12}],
    [0,  6, 16, x_{7}],
\\~
    [0,  9, 15, x_{9}],
    [0, 10, 18, x_{10}],
    [0, 12, 19, x_{11}],
    [0, 14, 28, 11],
    [0, 15, 32, x_{1}],
    [0, 16, 25, x_{8}],
\\~
    [0, 19,  1, x_{2}],
    [0, 21,  9, x_{4}],
    [0, 28,  5, x_{3}],
    [0, 31, 12, x_{5}],
    [0, 32, 17, x_{6}]
\end{array}\]

For $n=17$ and $u=16, 20$, and $22$, the
designs are provided in Appendix B1. \e

\begin{Lemma}\label{3n19u1}
  An HSD$(3^{19} u^1)$ exists for all feasible values of $u$.
\end{Lemma}
\p Lemma~\ref{necessary} tells that $u$ is odd and $0 \le u \le 27$.
Lemma~\ref{GDD3nu1}$(c)$ covers the cases when $u$ is odd and $u \equiv 3~(mod~6)$.
For $u=1, 5, 7$, $11$, and $13$, the designs are
covered by Lemmas~\ref{lem3n11}, \ref{lem3n51}, \ref{lem3n71}, \ref{lem3n111}, and
\ref{lem3n131},
respectively.  For $u=17, 19, 23$, and $25$, the
the designs are provided in Appendix B2. \e

\begin{Lemma}\label{3n23u1}
An HSD$(3^{23} u^1)$ exists if $u$ is odd and $17 \le u \le 23$.
\end{Lemma}
\p Lemma~\ref{GDD3nu1}$(c)$ covers the cases when $u$ is odd and
$u \equiv 0~(mod~3)$. For other values of $u$,
the designs are provided in Appendix B3.
\e

\section{Main Results}

Combining the results from the previous section, it is easy to
establish our main results.

\begin{Theorem}\label{3nu1to9}
  For $0 \le u \le 15$, an HSD$(3^n u^1)$ exists if and only if
  $n \ge 1+2u/3$, $n \ge 4$, and $n(n +2u-1)\equiv 0~(mod~4)$.
\end{Theorem}
\p The necessary condition comes from Lemma~\ref{necessary}.  For
$u = 0, 3$, the theorem is true because of Theorem~\ref{known0}$(d)$.
For $u = 1,2,4,5,7,8, 10, 11$, we have Lemmas~\ref{lem3n11}, \ref{lem3n21},
\ref{lem3n41}, \ref{lem3n51}, \ref{lem3n71}, \ref{lem3n81},
\ref{lem3n101}, \ref{lem3n111}, \ref{lem3n131}, and \ref{lem3n141}.
For $u = 6, 9, 12, 15$, Lemma~\ref{GDD3nu1} applies.  \e

\begin{Theorem}
  For $0 \le u \le n$, an HSD$(3^n u^1)$
  exists if and only if $n(n + 2u -1) \equiv 0~(mod~4)$ and $n \ge 4$,
  with possible exceptions of $n = 29, 43$.
\end{Theorem}
\p
By Lemma~\ref{necessary},
when $n \equiv 2~(mod~4)$, no HSD$(3^n u^1)$ exists.
For $4 \le n\le 15$, we have an HSD$(3^n u^1)$ for $0\le u \le n$ from
Theorem~\ref{3nu1to9}.
For $u\ge 16$ and $n \ge 12$, we consider $n = 4m+t$ for $m\ge 3$ and $t=0,1,3$.
\begin{itemize}
\item When $t=0$, an HSD$(3^n u^1)$ exists for $u \le 6m-2$
by Theorem~\ref{TD6m4m}.

\item When $t=1$, $u$ must be even by Lemma~\ref{necessary}.
Lemma~\ref{3n17u1} covers the cases when $n=13, 17$.
For $n\ge 21$, an HSD$(3^n u^1)$ exists for even $u \le n$ by
Lemma~\ref{TD6m4m1}, with possible exception of $n=29$.

\item When $t=3$, $u$ must be odd by Lemma~\ref{necessary}.
Lemmas~\ref{3n15u1}, \ref{3n19u1}, \ref{3n23u1}, and \ref{3n27u1} cover the cases
for $n \in \{ 15, 19, 23, 27 \}$, respectively.
For $n\ge 31$, an HSD$(3^n u^1)$ exists for odd $u \le n$ by
Lemma~\ref{TD6m4m3}, with possible exception of $n = 43$. \e
\end{itemize}

\ignore{
  Use HSD$(3^{32}s^1)$, $0\le s\le 46$, and HSD$(3^{11}u^1)$, $u$ odd and $1\le u \le 13$,
  we may obtain an HSD$(3{43}u^1)$ for $u$ odd and $1\le u \le 13$.
}

Besides the two possible exceptions in the above theorem, the cases of
HSD$(3^n u^1)$, where $max(15, n) < u \le (3n-3)/2$,
are left as a subject of further research.


Finally, we would like to report that six possible cases are removed
from Theorem~\ref{known0}$(f)$ and an updated version is given as
follows:

\begin{Theorem}
For $0 \le u \le 36$, an HSD$(4^nu^1)$ exists if
and only if $n \ge 4$ and
$0\le u \le 2n -2$,  with possible
exceptions of $(n, u) \in \{ (19, 29), (22, 33), (22, 35) \}$.
\end{Theorem}
\p
HSDs of types $(4^nu^1)$ for
$n = 19$ and $u \in \{ 30, 31, 33, 34, 35\}$, and
$(n, u) = (22, 34)$ are available in Appendix D.
\e

\subsection*{Acknowledgments}
The author wishes to thank Frank Bennett for early collaboration on
related projects which provided the foundation for the completion of
the work reported here.
The first author is supported by NSFC (No. 12161010).

\section*{Appendix}

The entire appendix which contains all the directed constructions referred in the paper
is available at
\begin{verbatim}
    https://homepage.divms.uiowa.edu/~hzhang/HSD3nu1Appendix.pdf
\end{verbatim}

\newpage

\section*{Appendix}

Here we list some HSDs which are used in the previous sections and
obtained by computer. The point set of an HSD$(h^nu^1)$ consists of
$Z_{hn}$ and $u$ infinite points which are denoted by alphabet.  For
simplicity, we only list the starter blocks.
We also use the $+k$ method to develop blocks, which means
that we add $k$ (mod $hn$) to each point of the starter blocks to
obtain all blocks. Starter blocks of short orbit are in bold cases
(see Example~\ref{ex3821}).

\subsubsection*{A1 HSD$(3^n 1^1)$}

{\footnotesize\[\begin{array}{l}
h = 3, n = 8, u = 1 ~(+4~ mod~ 24):
\\~[0,  1,  2, 19], [0,  2,  6,  9], [0,  3, 21, x_1], [0,  4,  5,  7], [0,  5,  1, 11], [0,  6, 15, 20], [0,  7, 18,  3],
\\~[0,  9, 22, 23], [0, 10, 12, 22], [0, 11, 23, 12], [0, 12, 17,  5], [0, 14,  7, 10], [0, 15, 11, 13],
\\~[0, 17, 14, 21],
[0, 18,  9, 15], [0, 21,  4,  2], [0, 23, 13, 18], [1,  2, 23, 19], [1,  5, 10, 15], [1, 10, 22, 13],
\\~[1, 12,  2, x_1],
[1, 15, 13,  3], [1, 19, 14, 18], [2, 13, 11, x_1], [2, 14, 15,  3], [3, 14,  8, x_1]
\\\\
h = 3, n = 11, u = 1 ~(+1~ mod~ 33):
\\~[0,  1,  2, x_1], [0,  2, 21, 29], [0,  3, 28, 16], [0,  4, 30,  3], [0,  5, 25, 15], [0,  7,  9, 24], [0,  9, 14, 27],
[0, 14,  4, 21]
\\
\\h = 3, n = 12, u = 1 ~(+2~ mod~ 36):
\\~
{\bf [0, 18, 13, 31]},
{\bf [0, 25,  7, 18]},
{\bf [0, 33, 18, 15]},
[0,  1, 21, 26],
[0,  2, 30,  3],
[0,  3, 16, 30],
[0,  4,  2,  1],
\\~
[0,  5,  3, 22],
[0,  6,  5, 21],
[0,  7, 26, 11],
[0,  8, 27, 16],
[0, 10, 35,  5],
[0, 13,  4, 27],
[0, 15, 22, 13],
\\~
[0, 16, 23, 33],
[0, 19, 29, x_1],
[1,  3, 31, 23],
[1,  5, 11, 33],
[1,  8, 12, x_1]
\\
\\h = 3, n = 15, u = 1 ~(+1~ mod~ 45):
\\~[0,  1, 23, 41], [0,  3,  1,  9], [0,  4, 21, 11], [0,  5, 37, 25], [0,  6, 11, 35], [0,  7, 31,  3], [0,  9, 35, 12],
\\~[0, 11, 17, 37], [0, 13, 12, 31], [0, 14, 43, 27], [0, 43, 36, x_1]
\end{array}\]}

\subsubsection*{A2 HSD$(3^n 2^1)$}

{\footnotesize\[\begin{array}{l}
h = 3, n = 9, u = 2 ~(+1~ mod~ 27):
\\~[0,  1, 22, x_1], [0,  2,  6, 17], [0,  3,  2,  7], [0,  4,  7, 14], [0,  6, 11, 19], [0, 12, 24, 11], [0, 17, 19, x_2]
\\\\
h = 3, n = 12, u = 2 ~(+2~ mod~ 36):
\\~
{\bf [0,  5, 23, 18]},
{\bf [0, 11, 18, 29]},
{\bf [0, 18, 27,  9]},
[0,  2, 22,  6],
[0,  3,  7, 14],
[0,  4,  8, 35],
[0,  6, 13,  3],
\\~
[0,  7, 30, 15],
[0,  8, 29, x_2],
[0,  9,  1, x_1],
[0, 10, 21, 26],
[0, 13, 19, 27],
[0, 14, 15, 17],
[0, 15, 17, 31],
\\~
[0, 17, 34, 23],
[0, 19,  5,  8],
[0, 23, 26, 25],
[1,  0,  2, x_1],
[1,  7, 27, 11],
[1, 33, 28, x_2]
\\\\
h = 3, n = 13, u = 2 ~(+1~ mod~ 39):
\\~[0,  1, 10, 16], [0,  2,  1, 10], [0,  3, 23, 12], [0,  4,  2,  7], [0,  7, 28, x_2], [0,  8, 22,  3], [0, 12, 20,  6],
[0, 15, 32, 11],
\\~[0, 16, 12, 34],
[0, 29, 14, x_1]
\\\\
h = 3, n = 17, u = 2 ~(+1~ mod~ 51):
\\~[0,  1, 29, 32], [0,  2, 39, x_1], [0,  4, 42, 27], [0,  5, 16,  3], [0,  6, 36, 12],
[0,  7,  1, 10], [0,  8, 37, 16], [0, 10, 30,  5],
\\~
[0, 11, 18, 49], [0, 12, 11, x_2], [0, 14, 10, 33], [0, 16, 25, 43], [0, 19,  4, 26]
\end{array}\]}

\subsubsection*{A3 HSD$(3^n 4^1)$}

{\footnotesize\[\begin{array}{l}
h = 3, n = 8, u = 4 ~(+2~ mod~24):
\\~
{\bf [0,  9, 21, 12]},
{\bf [0, 11, 12, 23]},
{\bf [0, 12,  7, 19]},
[0,  3, 17, x_1],
[0,  4, 11, 18],
[0,  5,  1, 15],
[0,  6, 19, x_2],
\\~
[0,  7, 13, 22],
[0, 10, 15,  9],
[0, 13,  2, x_4],
[0, 21,  4,  3],
[0, 22, 18, x_3],
[1,  0, 21, x_4],
[1,  6, 20, x_1],
\\~[1, 21, 23, x_3], [1, 23, 22, x_2]
\\
\\h = 3, n = 9, u = 4 ~(+1 ~mod~ 27):
\\~[0,  1,  2, x_4], [0,  3, 11, 16], [0,  6, 26, 14], [0,  7,  4, 12], [0, 10, 12, x_2], [0, 11, 17,  4],
 [0, 23,  6, x_1],
\\~[0, 25,  3, x_3]
\\
\\h = 3, n = 13, u = 4 ~(+1 ~mod~ 39):
\\~[0,  1,  2, x_4], [0,  3, 25, 11], [0,  4, 22, 27], [0,  6, 10, x_1], [0,  7,  9, x_2], [0,  8, 28,  9],
[0,  9, 34, 24],
\\~[0, 11, 19,  7], [0, 15, 21,  3], [0, 16, 32, 10], [0, 37,  3, x_3]
\\
\\h = 3, n = 17, u = 4 ~(+1~ mod~ 51):
\\~[0,  1,  7, 44], [0,  2, 28, 47], [0,  3, 46, 41], [0,  4, 40, 18], [0,  6, 33,  9], [0,  7,  9, 39], [0,  8, 30, x_3],
\\~[0,  9, 25, 13], [0, 10, 50, 32], [0, 11, 31, x_4], [0, 13, 16, x_1], [0, 15, 36,  5], [0, 16,  2, 28], [0, 28, 27, x_2]
\end{array}\]}

\subsubsection*{A4 HSD$(3^n 5^1)$}

{\footnotesize\[\begin{array}{l}
h = 3, n = 7, u = 5  ~(+1 ~mod~ 21):
\\~[0,  1,  2, x_5], [0,  3,  6, 12], [0,  5, 20, 10], [0,  8, 10, x_3], [0, 12,  8, x_1], [0, 17,  4, x_2],
 [0, 19,  3, x_4]
\\
\\h = 3, n = 8, u = 5 ~(+2~ mod~ 24):
\\~[0,  1, 20, 22], {\bf [0,  3, 15, 12]}, [0,  4, 18, x3], [0,  9,  2, 15], [0, 10, 19, 17], {\bf [0, 11, 12, 23], [0, 12, 23, 11]},
\\~[0, 15,  5, x1], [0, 17, 11, x2], [0, 18, 21, x4], [0, 19, 17, 21], [0, 23, 10, x5], [1,  4,  5, x5], [1, 11, 16, x4],
\\~[1, 18, 12, x2], [1, 19, 15, x3], [1, 20,  0, x1]
\\
\\h = 3, n = 11, u = 5 ~(+1 ~mod~ 33):
\\~[0,  1,  2, x_5], [0,  3, 23, 15], [0,  4, 28, 21], [0,  5, 13, x_1], [0,  6,  8, x_3], [0,  9,  6, 27],
\\~[0, 10, 14, x_2], [0, 13, 32, 17], [0, 14, 24,  7], [0, 31,  3, x_4]
\\\\
h = 3, n = 19, u = 5 ~(+1~ mod~ 57):
\\~[0,  1, 54, 37], [0,  2, 26, x4], [0,  3, 30, 43], [0,  4,  7, 16], [0,  5, 53, 11], [0,  6, 32, 50], [0,  7, 47,  6],
\\~[0,  8, 18, 49], [0, 10, 23, 52], [0, 11,  9, 39], [0, 12, 35, x1], [0, 20, 49, 15], [0, 21, 20, 45], [0, 22, 43, x2],
\\~[0, 33, 55, x3], [0, 43, 11, x5]
\\
\\h = 3, n = 23, u = 5  ~(+1~ mod~ 69):
\\~[0,  1,  2, x5], [0,  2, 64, 37], [0,  3, 62, 57], [0,  4,  6, x4], [0,  6, 66, 53], [0,  7, 68, 31], [0,  8, 26, 50],
\\~[0,  9, 12, x3], [0, 10, 30, 42], [0, 11, 52, 22], [0, 14, 38, 55], [0, 15, 51, 25], [0, 16, 20, x2], [0, 18, 61, 39],
\\~[0, 19, 40,  6], [0, 20, 55, 11], [0, 21, 50, 17], [0, 28, 33, x1], [0, 29, 44, 13]
\end{array}\]}

\subsubsection*{A5 HSD$(3^n 7^1)$}

{\footnotesize\[\begin{array}{l}
h = 3, n = 7, u = 7 ~(+1 ~mod~ 21):
\\~[0,  1,  2,x_7], [0, 2, 4,x_6], [0, 3,18  9], [0, 4, 9,x_4], [0, 5,  8 x_2],
[0,  6, 10, x_5], [0,  8, 16, x_3], [0, 10, 20, x_1]
\end{array}\]}

{\footnotesize\[\begin{array}{l}
h = 3, n = 8, u = 7 ~(+4 ~mod~ 24):
\\~[0,  2,  6, x_7],
                  [0,  6,  3, 17],
                  [0,  7,  9, x_6],
                  [0, 10,  1, 19],
                  [0, 12, 11, 23],
                  [0, 13, 12,  1],
                  [0, 14,  2, 12],
\\~[0, 15, 14, x_5],
                  [0, 17, 19, 14],
                  [0, 18,  5, x_4],
                  [0, 19,  4, x_2],
                  [0, 20,  7, x_3],
                  [0, 22, 15, x_1],
                  [0, 23, 10,  5],
\\~[1,  0,  4, x_7],
                  [1,  2,  7, x_6],
                  [1,  4, 11, x_4],
                  [1, 10, 12, x_3],
                  [1, 13, 10, 22],
                  [1, 15,  3, 13],
                  [1, 16, 18, x_2],
\\~[1, 20,  5, x_5],
                  [1, 21,  8, x_1],
                  [2,  3,  7, x_7],
                  [2,  5,  8, x_6],
                  [2,  9, 15, x_2],
                  [2, 11, 21, x_1],
                  [2, 13,  6, x_4],
\\~[2, 19,  1, x_3],
                  [2, 22,  4, x_5],
                  [2, 23, 14, 11],
                  [3,  0,  6, x_6],
                  [3,  1,  5, x_7],
                  [3,  5, 10, x_3],
                  [3,  7,  4, x_4],
\\~[3, 14, 13, x_2],
              [3, 16,  2, x_1],
                  [3, 21,  7, x_5]
                  \\
\\h = 3, n = 11, u = 7 ~(+1 ~mod~ 33):
\\~[0,  2, 30, x_7], [0,  3, 10, 30], [0,  4, 25, x_5], [0,  5, 21, x_6], [0,  8, 18,  9], [0, 10, 29, 15],
 [0, 12, 14, 29],
\\~[0, 16, 24, x_4, [0, 26, 20, x_1], [0, 27,  7, x_2], [0, 32, 31, x_3]
\\
\\h = 3, n = 19, u = 7 ~(+1 ~mod~ 57):
\\~[0,  1,  2, x_7],
                  [0,  2,  4, x_6],
                  [0,  3, 50, 37],
                  [0,  4, 25, 35],
                  [0,  5, 16,  8],
                  [0,  6, 48, 28],
                  [0,  7, 10, x_5],
\\~[0,  9, 21, 39],
                  [0, 11, 17, x_2],
                  [0, 12, 52, 25],
                  [0, 14, 18, x_4],
                  [0, 15, 42, 26],
                  [0, 17, 51, 16],
                  [0, 21, 45, 13],
\\~[0, 23, 37,  9],
                  [0, 24, 29, x_3],
                     [0, 26, 33, x_1]
\\
\\h = 3, n = 23, u = 7 ~(+1 ~mod~ 69):
\\~[0,  1, 41,  7],
                  [0,  2, 10,  5],
                  [0,  3, 42, 22],
                  [0,  6,  7, 38],
                  [0,  7, 64, 25],
                  [0,  8, 67, 51],
                  [0,  9, 35, 63],
\\~[0, 10, 12, 24],
                  [0, 11, 25, 49],
                  [0, 13, 16, 33],
                  [0, 14, 48, x_5],
                  [0, 15, 68, 28],
                  [0, 18, 60, x_6],
                  [0, 19, 30, 52],
\\~[0, 21, 58, x_1],
                  [0, 26, 17, 50],
                  [0, 32, 47, x_4],
                  [0, 42, 29, x_2],
                  [0, 44, 65, x_3],
                   [0, 65, 61, x_7]
\end{array}\]}

\subsubsection*{A6 HSD$(3^n 8^1)$}

{\footnotesize\[\begin{array}{l}
h = 3, n = 8, u = 8 ~(+4 ~mod~ 24):
\\~[0,  2,  6, x_8],
                  [0,  7,  9, x_7],
                  [0,  9,  4, x_2],
                  [0, 10, 19, x_4],
                  [0, 12, 10, 22],
                  [0, 13, 12,  1],
                  [0, 14,  3, 13],
\\~[0, 15, 14, x_6],
                  [0, 17, 11,  2],
                  [0, 18,  5, x_5],
                  [0, 19,  7, 12],
                  [0, 20, 15, x_1],
                  [0, 22,  1, x_3],
                  [1,  0,  4, x_8],
\\~[1,  2,  7, x_7],
                  [1,  4, 11, x_5],
                  [1,  6,  3, x_2],
                  [1, 10, 22, 13],
                  [1, 11, 12, x_3],
                  [1, 13, 15,  3],
                  [1, 19,  8, x_4],
\\~[1, 20,  5, x_6],
                  [1, 21,  6, x_1],
                  [2,  3,  7, x_8],
                  [2,  5,  8, x_7],
                  [2,  7, 14, 19],
                  [2,  9, 15, x_3],
                  [2, 13,  6, x_5],
\\~[2, 15,  0, x_1],
                  [2, 20, 13, x_4],
                  [2, 22,  4, x_6],
                  [2, 23,  9, x_2],
                  [3,  0,  6, x_7],
                  [3,  1,  5, x_8],
                  [3,  4, 18, x_3],
\\~[3,  5, 10, x_4],
                  [3,  7,  4, x_5],
                  [3, 10,  9, x_1],
                  [3, 16,  2, x_2],
                  [3, 21,  7, x_6]
\\\\
h = 3, n = 9, u = 8 ~(+1 ~mod~ 27):
\\~[0,  1,  2, x_8], [0,  2,  4, x_7], [0,  5, 16, x_1], [0,  6, 20,  7], [0,  8, 12, x_4], [0, 10, 13, x_5],
 [0, 11, 17,  5],
\\~[0, 20,  3, x_3], [0, 23,  8, x_2], [0, 24,  5, x_]
\\\\
h = 3, n = 13, u = 8 ~(+1 ~mod~ 39):
\\~[0,  1, 19, 17], [0,  3, 17,  7], [0,  4,  5, 37], [0,  5, 15, x_8], [0,  6,  2, x_2], [0,  9, 14, 30], [0, 31,  7, x_7],
[0, 11,  8, 20],
\\~[0, 15, 38, x_6], [0, 17, 29, x_3], [0, 19, 11, x_1], [0, 21, 27, x_4], [0, 25, 36, x_5]
\\\\
h = 3, n = 17, u = 8 ~(+1 ~mod~ 51):
\\~[0,  3, 35, 23], [0,  4, 29, x_8], [0,  5, 26, 15], [0,  7, 11, 39], [0,  9, 33,  3], [0, 10, 43,  7],  [0, 43,  5, x_3],
\\~[0, 13, 44, x_2], [0, 14, 30, 50], [0, 16, 38, 14], [0, 18,  9, 41], [0, 26, 24, x_1], [0, 29, 37, x_4],
\\~[0, 45, 50, x_7], [0, 49,  4, x_5], [0, 50, 39, x_6]
\\\\
h = 3, n = 29, u = 8 ~(+1 ~mod~ 87):
\\~
    [0,  1,  2, x_8],
    [0,  2,  4, x_7],
    [0,  3, 79, 51],
    [0,  4, 20, 60],
    [0,  5, 14, x_1],
    [0,  6, 65, 50],
    [0,  7, 10, x_6],
\\~
    [0,  8, 38, 70],
    [0,  9, 51, 33],
    [0, 10, 66, 25],
    [0, 11, 52, 68],
    [0, 12, 50, 72],
    [0, 13, 47, 66],
    [0, 14, 18, x_5],
\\~
    [0, 17, 81, 14],
    [0, 21, 82, 20],
    [0, 23, 28, x_4],
    [0, 24, 64, 13],
    [0, 26, 19, 52],
    [0, 27, 75, 18],
    [0, 31, 74, 24],
\\~
    [0, 34, 42, x_2],
    [0, 35, 41, x_3],
    [0, 38, 16, 55],
    [0, 42, 32, 75]
\end{array}\]}

\subsubsection*{A7 HSD$(3^n 10^1)$}

{\footnotesize\[\begin{array}{l}
h = 3, n = 8, u = 10 ~(+8 ~mod~ 24):
\\~
   [0,  1, 11, x_1],
   [0,  3,  6, x_2],
   [0,  6, 15, x_{10}],
   [0,  7, 14, x_9],
   [0, 10, 23, x_5],
   [0, 11,  5, x_7],
   [0, 14, 21, x_8],
\\~
   [0, 15, 18, x_6],
   [0, 18, 22, 11],
   [0, 19, 17, x_4],
   [0, 23,  1, x_3],
   [1,  2, 16, x_2],
   [1,  4, 14, x_7],
   [1,  5, 15, x_9],
\\~
   [1,  8, 22, x_{10}],
   [1, 10, 13, x_3],
   [1, 11, 23, x_6],
   [1, 12, 21, x_5],
   [1, 13,  2, x_1],
   [1, 15,  4, x_4],
   [1, 16,  7, x_8],
\\~
   [1, 21, 20,  7],
   [2,  0, 20, x_7],
   [2,  3,  9, x_9],
   [2,  4, 23, x_1],
   [2,  5,  6, x_4],
   [2,  7,  5, x_2],
   [2,  9, 16, x_{10}],
\\~
   [2, 11, 22, x_5],
   [2, 12, 17, x_6],
   [2, 13,  0, x_8],
   [2, 14, 15, x_3],
   [3,  1, 14, x_8],
   [3,  4,  0, x_4],
   [3,  5,  9, x_7],
\\~
   [3,  6,  1, x_2],
   [3,  9,  6, x_3],
   [3, 10, 17, x_{10}],
   [3, 12, 18, x_9],
   [3, 13, 22, x_6],
   [3, 15,  8, x_1],
   [3, 21,  2, x_5],
\\~
   [4,  0, 22, x_1],
   [4,  5, 19, x_2],
   [4,  7, 11, x_8],
   [4,  8,  3, x_4],
   [4,  9,  0, x_6],
   [4, 10, 16, x_9],
   [4, 11, 18, x_{10}],
\\~
   [4, 13,  8, x_3],
   [4, 16, 17, x_5],
   [4, 22,  2, x_7],
   [5,  0, 12, x_9],
   [5,  6,  4, x_1],
   [5,  7,  3, x_5],
   [5,  8, 11, x_3],
\\~
   [5, 10,  1, x_8],
   [5, 12, 19, x_{10}],
   [5, 14, 20, x_6],
   [5, 15, 16, x_7],
   [5, 16, 18, x_2],
   [5, 22, 10, x_4],
   [6,  1, 13, x_9],
\\~
   [6,  2,  1, x_1],
   [6,  4,  7, x_2],
   [6,  8,  5, x_6],
   [6,  9, 23, x_7],
   [6, 10, 15, x_4],
   [6, 11, 20, x_3],
   [6, 15, 21, x_{10}],
\\~
   [6, 17,  4, x_5],
   [6, 20, 18, x_8],
   [7,  1,  4, x_2],
   [7,  3,  5, x_1],
   [7,  6,  0, x_5],
   [7,  9, 13, x_4],
   [7, 10, 11, x_7],
\\~
   [7, 11, 12, x_8],
   [7, 12,  2, x_3],
   [7, 13, 20, x_{10}],
   [7, 14, 19, x_9],
   [7, 18,  3, x_6]
\\\\
h = 3, n = 9, u = 10 ~(+1 ~mod~ 27):
\\~
   [0,  1,  2, x_{10}],
   [0,  2,  4, x_9],
   [0,  3,  6, x_8],
   [0,  4, 12, x_4],
   [0,  5, 20, 13],
   [0,  6, 10, x_5],
   [0,  8, 13, x_7],
\\~
   [0, 10, 16, x_6],
   [0, 11,  1, x_3],
   [0, 12,  5, x_2],
   [0, 13, 24, x_1]
\\
\\h = 3, n = 13, u = 10 ~(+1 ~mod~ 39):
\\~
   [0,  1,  2, x_{10}],
   [0,  2,  4, x_9],
   [0,  3,  6, x_8],
   [0,  4, 23, x_2],
   [0,  5, 12, x_3],
   [0,  6, 36, 11],
   [0,  8, 31, 22],
\\~
   [0, 10, 15, x_5],
   [0, 11, 19, x_4],
   [0, 12, 18, x_6],
   [0, 15,  1, 24],
   [0, 17, 29, 10],
   [0, 18, 22, x_7],
   [0, 32, 11, x_1]
\\
\\h = 3, n = 17, u = 10 ~(+1 ~mod~ 51):
\\~
   [0,  1,  2, x_{10}],
   [0,  2,  4, x_9],
   [0,  3,  6, x_8],
   [0,  5, 48, 24],
   [0,  6, 16, x_1],
   [0,  7, 39, 28],
   [0,  8, 38, 26],
\\~
   [0,  9, 44, 29],
   [0, 10, 14, x_7],
   [0, 13, 46, 14],
   [0, 14, 23, x_2],
   [0, 16, 22, x_4],
   [0, 18, 25, x_3],
   [0, 20, 40, 12],
\\~
   [0, 21, 10, 36],
   [0, 22, 27, x_6],
   [0, 47,  9, x_5]
\\\\
h = 3, n = 29, u = 10 ~(+1 ~mod~ 87):
\\~
    [0,  4, 53, 16],
    [0,  5, 64, 23],
    [0,  6, 62, 20],
    [0,  7, 55, 24],
    [0,  8, 84,  7],
    [0, 11, 56, 68],
    [0, 14, 80,  9],
\\~
    [0, 15, 48, 75],
    [0, 17, 63, 83],
    [0, 18, 35, 65],
    [0, 19, 11, 32],
    [0, 22, 65,  6],
    [0, 23, 50, x_7],
    [0, 24, 19, 52],
\\~
    [0, 32, 45, x_6],
    [0, 35, 72, 25],
    [0, 36, 26, 69],
    [0, 38, 78, 30],
    [0, 53, 67, x_4],
    [0, 61, 46, x_1],
    [0, 62, 36, x_5],
\\~
    [0, 74, 38, x_2],
    [0, 78, 44, x_3],
    [0, 84, 81, x_8],
    [0, 85, 83, x_9],
    [0, 86, 85, x_{10}]
\end{array}\]}

\subsubsection*{A8 HSD$(3^n 11^1)$}

{\footnotesize\[\begin{array}{l}
h = 3, n = 11, u = 11 ~(+1 ~mod~ 33):
\\~
   [0,  1,  2, x_{11}],
   [0,  2,  4, x_{10}],
   [0,  3,  6, x_9],
   [0,  4, 13, 20],
   [0,  8, 14, x_6],
   [0,  9, 32, 18],
   [0, 10, 15, x_7],
\\~
   [0, 13, 30, x_2],
   [0, 17, 21, x_5],
   [0, 18, 10, x_3],
   [0, 21, 28, x_1],
   [0, 27,  8, x_4],
   [0, 28,  7, x_8]
\\\\
h = 3, n = 19, u = 11 ~(+1 ~mod~ 57):
\\~
   [0,  1,  2, x_{11}],
   [0,  2,  4, x_{10}],
   [0,  3,  6, x_9],
   [0,  4, 16, x_1],
   [0,  5, 11, x_5],
   [0,  6, 15, x_3],
   [0,  7, 33, 17],
\\~
   [0,  8, 29, 20],
   [0, 10, 14, x_8],
   [0, 11, 44, 29],
   [0, 12, 23, 44],
   [0, 13, 52, 34],
   [0, 14, 22, x_4],
   [0, 17, 49, 14],
\\~
   [0, 20, 27, x_6],
   [0, 23, 50, 22],
   [0, 24,  9, 41],
   [0, 26, 31, x_7],
   [0, 27, 37, x_2]
\end{array}\]}

{\footnotesize\[\begin{array}{l}
h = 3, n = 23, u = 11 ~(+1 ~mod~ 69):
\\~
   [0,  1,  2, x_{11}],
   [0,  2,  4, x_{10}],
   [0,  3,  6, x_9],
   [0,  4, 43, 49],
   [0,  5, 29, 42],
   [0,  7, 44, 14],
   [0,  8, 61, 35],
\\~
   [0,  9, 16, x_6],
   [0, 10, 14, x_8],
   [0, 11, 58, 26],
   [0, 12, 60, 31],
   [0, 14, 22, x_2],
   [0, 15, 48, 28],
   [0, 16, 66, 15],
\\~
   [0, 17, 57, 12],
   [0, 19, 30, x_3],
   [0, 21, 31, x_4],
   [0, 22, 28, x_5],
   [0, 25, 34, x_1],
   [0, 27, 10, 44],
   [0, 28, 33, x_7],
\\~
   [0, 31, 49, 13]
\end{array}\]}

\subsubsection*{A9 HSD$(3^n 13^1)$}

{\footnotesize\[\begin{array}{l}
h = 3, n = 11, u = 13 ~(+1~ mod~ 33):
\\~
    [0,  1,  2, x_{13}],
    [0,  2,  4, x_{12}],
    [0,  3,  6, x_{11}],
    [0,  4,  8, x_{10}],
    [0,  5, 13, x_{7}],
    [0,  6, 23, x_{4}],
    [0,  7, 28, x_{2}],
\\~
    [0,  8, 15, x_{8}],
    [0,  9, 14, x_{6}],
    [0, 10, 16, x_{9}],
    [0, 12, 21, x_{5}],
    [0, 15, 30, 14],
    [0, 19,  9, x_{3}],
    [0, 20,  7, x_{1}]
\\\\
h = 3, n = 19, u = 13 ~(+1~ mod~ 57):
\\~
    [0,  1,  2, x_{13}],
    [0,  2,  4, x_{12}],
    [0,  3,  6, x_{11}],
    [0,  4,  8, x_{10}],
    [0,  6, 17, x_3],
    [0,  7, 25, 36],
    [0,  8, 18, x_4],
\\~
    [0,  9, 16, x_5],
    [0, 10, 52, 32],
    [0, 12, 56, 26],
    [0, 14, 20, x_9],
    [0, 15, 48, 30],
    [0, 16, 44, 23],
    [0, 17, 26, x_2],
\\~
    [0, 22, 10, 43],
    [0, 23, 45, 20],
    [0, 26, 34, x_6],
    [0, 28, 33, x_8],
    [0, 44, 27, x_1],
    [0, 52, 11, x_7]
\\\\
h = 3, n = 23, u = 13 ~(+1~ mod~ 69):
\\~
    [0,  1,  2, x_{13}],
    [0,  2,  4, x_{12}],
    [0,  3,  6, x_{11}],
    [0,  4,  8, x_{10}],
    [0,  7, 59, 47],
    [0,  8, 62, 33],
    [0,  9, 53, 35],
\\~
    [0, 10, 50, 13],
    [0, 11, 52, 32],
    [0, 13, 18, x_9],
    [0, 14, 28, 45],
    [0, 15, 27, x_1],
    [0, 16, 24, x_7],
    [0, 19, 58, 31],
\\~
    [0, 21, 30, x_2],
    [0, 22, 64, 21],
    [0, 24, 34, x_3],
    [0, 25, 32, x_6],
    [0, 28,  9, 43],
    [0, 30, 36, x_8],
    [0, 31, 47, 11],
\\~
    [0, 63, 14, x_4],
    [0, 64, 13, x_5]
\end{array}\]}

\subsubsection*{A10 HSD$(3^n 14^1)$}

{\footnotesize\[\begin{array}{l}
h = 3, n = 13, u = 14 ~(+1~ mod~ 39):
\\~
    [0,  1,  2, x_{14}],
    [0,  2,  4, x_{13}],
    [0,  3,  6, x_{12}],
    [0,  4,  8, x_{11}],
    [0,  5, 28, x_{4}],
    [0,  6, 16, x_{6}],
    [0,  7, 15, x_{7}],
\\~
    [0,  8, 29, 20],
    [0, 10, 17, x_{8}],
    [0, 11, 38, x_{3}],
    [0, 12, 18, x_{10}],
    [0, 14, 19, x_{9}],
    [0, 15, 30, x_{2}],
    [0, 16, 25, x_{5}],
\\~
    [0, 17, 34, 14],
    [0, 18,  7, x_{1}]
\\\\
h = 3, n = 17, u = 14 ~(+1~ mod~ 51):
\\~
    [0,  1,  2, x_{14}],
    [0,  2,  4, x_{13}],
    [0,  3,  6, x_{12}],
    [0,  4,  8, x_{11}],
    [0,  6, 16, x_4],
    [0,  7, 40, 20],
    [0,  8, 46, 27],
\\~
    [0,  9, 44, 23],
    [0, 10, 21, x_3],
    [0, 11, 33, x_2],
    [0, 12, 20, x_5],
    [0, 13, 50, 25],
    [0, 14, 23, x_6],
    [0, 15,  3, 30],
\\~
    [0, 16, 22, x_8],
    [0, 18, 25, x_7],
    [0, 22, 27, x_{10}],
    [0, 23, 42, x_1],
    [0, 46, 10, x_9]
\\\\
h = 3, n = 29, u = 14 ~(+1~ mod~ 87):
\\~
   [0,  1,  2, x_{14}],
   [0,  2,  4, x_{13}],
   [0,  3,  6, x_{12}],
   [0,  4,  8, x_{11}],
   [0,  5, 76, 57],
   [0,  6, 80, 53],
   [0,  7, 77, 56],
\\~
   [0,  8, 70, 54],
   [0,  9, 27, 37],
   [0, 11, 62, 47],
   [0, 12, 20, x_8],
   [0, 13, 18, x_{10}],
   [0, 14, 78, 46],
   [0, 17, 31, x_3],
\\~
   [0, 18, 33, x_2],
   [0, 20, 30, x_5],
   [0, 22, 34, x_1],
   [0, 23, 75, 44],
   [0, 24, 84, 38],
   [0, 25, 86, 48],
   [0, 26, 68, 21],
\\~
   [0, 28, 39, x_6],
   [0, 30, 37, x_7],
   [0, 33, 65, 20],
   [0, 34, 43, x_4],
   [0, 35, 63, 19],
   [0, 36, 42, x_9],
   [0, 37, 61, 22]

\end{array}\]}

\subsubsection*{Appendix B1 HSD$(3^{17} u^1)$}

{\footnotesize\[\begin{array}{l}
h = 3, n = 17, u = 16 ~(+1~ mod~ 51):
\\~
    [0,  1,  2, x_{16}],
    [0,  2,  4, x_{15}],
    [0,  3,  6, x_{14}],
    [0,  4,  8, x_{13}],
    [0,  5, 10, x_{12}],
    [0,  6, 35, x_3],
    [0,  8, 21, x_4],
\\~
    [0,  9, 19, x_5],
    [0, 10, 33, x_2],
    [0, 11, 22, x_6],
    [0, 12, 20, x_7],
    [0, 13, 37, 12],
    [0, 14, 44, 26],
    [0, 15, 24, x_8],
\\~
    [0, 16, 23, x_{10}],
    [0, 19, 25, x_9],
    [0, 20,  5, 35],
    [0, 22, 42, 19],
    [0, 24, 38, x_1],
    [0, 44, 11, x_{11}]
\\\\
h = 3, n = 17, u = 20 ~(+1~ mod~ 51):
\\~
    [0,  1,  2, x_{20}],
    [0,  2,  4, x_{19}],
    [0,  3,  6, x_{18}],
    [0,  4,  8, x_{17}],
    [0,  5, 10, x_{16}],
    [0,  6, 12, x_{15}],
    [0,  8, 22, x_{8}],
\\~
    [0,  9, 20, x_{10}],
    [0, 10, 19, x_{11}],
    [0, 11, 21, x_{9}],
    [0, 15, 38, x_{3}],
    [0, 16, 24, x_{14}],
    [0, 18, 25, x_{13}],
    [0, 19, 44, 18],
\\~
    [0, 20,  5, x_{2}],
    [0, 21, 33, x_{5}],
    [0, 22, 42, 19],
    [0, 24, 11, x_{1}],
    [0, 37, 16, x_{6}],
    [0, 38, 14, x_{4}],
    [0, 39, 23, x_{7}],
\\~
    [0, 44, 15, x_{12}]
\\\\
h = 3, n = 17, u = 22 ~(+1~ mod~ 51):
\\~
    [0,  1,  2, x_{22}],
    [0,  2,  4, x_{21}],
    [0,  3,  6, x_{20}],
    [0,  4,  8, x_{19}],
    [0,  5, 10, x_{18}],
    [0,  6, 12, x_{17}],
    [0,  7, 14, x_{16}],
\\~
    [0, 10, 23, x_{9}],
    [0, 11, 22, x_{13}],
    [0, 13, 21, x_{11}],
    [0, 14, 24, x_{14}],
    [0, 16, 25, x_{15}],
    [0, 18, 38, x_{5}],
    [0, 19, 44, 18],
\\~
    [0, 20,  5, x_{4}],
    [0, 21, 33, x_{8}],
    [0, 22,  3, x_{2}],
    [0, 24, 40, x_{3}],
    [0, 28, 42, x_{1}],
    [0, 36, 15, x_{6}],
    [0, 39, 16, x_{7}],
\\~
    [0, 42, 20, x_{10}],
    [0, 43, 19, x_{12}]
\end{array}\]}

\subsubsection*{B2 HSD$(3^{19} u^1)$}

{\footnotesize\[\begin{array}{l}
h = 3, n = 19, u = 17 ~(+1~ mod~ 57):
\\~
    [0,  1,  2, x_{17}],
    [0,  2,  4, x_{16}],
    [0,  3,  6, x_{15}],
    [0,  4,  8, x_{14}],
    [0,  5, 10, x_{13}],
    [0,  6, 18, x_6],
    [0,  7, 22, x_4],
    [0,  8, 21, x_5],
\\~
    [0, 10, 20, x_7],
    [0, 11, 50, 26],
    [0, 12, 41, x_2],
    [0, 13, 45, 30],
    [0, 14, 23, x_8],
    [0, 16, 27, x_9],
    [0, 17, 54, 22],
\\~
    [0, 18, 26, x_{11}],
    [0, 20, 56, 34],
    [0, 21, 28, x_{10}],
    [0, 23,  9, x_1],
    [0, 26, 32, x_{12}],
    [0, 27, 44, 16],
    [0, 48, 24, x_3]
\\
\\h = 3, n = 19, u = 19 ~(+1~ mod~ 57):
\\~
    [0,  1,  2, x_{19}],
    [0,  2,  4, x_{18}],
    [0,  3,  6, x_{17}],
    [0,  4,  8, x_{16}],
    [0,  5, 10, x_{15}],
    [0,  6, 12, x_{14}],
    [0,  7, 18, x_9],
\\~
    [0,  9, 23, x_7],
    [0, 10, 22, x_8],
    [0, 12, 25, x_6],
    [0, 13, 52, 24],
    [0, 14, 50, 27],
    [0, 15, 24, x_{10}],
    [0, 16, 56, 32],
\\~
    [0, 17, 27, x_{11}],
    [0, 18, 26, x_{13}],
    [0, 21, 28, x_{12}],
    [0, 22, 37, x_4],
    [0, 25, 48, 22],
    [0, 27, 43, x_1],
    [0, 37, 17, x_2],
\\~
    [0, 46, 15, x_3],
    [0, 49, 21, x_5]
\\\\
h = 3, n = 19, u = 23 ~(+1~ mod~ 57):
\\~
    [0,  1,  2, x_{23}],
    [0,  2,  4, x_{22}],
    [0,  3,  6, x_{21}],
    [0,  4,  8, x_{20}],
    [0,  5, 10, x_{19}],
    [0,  6, 12, x_{18}],
    [0,  7, 14, x_{17}],
\\~
    [0,  9, 24, x_{10}],
    [0, 10, 22, x_{12}],
    [0, 11, 21, x_{13}],
    [0, 12, 23, x_{11}],
    [0, 13, 26, x_{8}],
    [0, 17, 46, x_{3}],
    [0, 18, 27, x_{16}],
\\~
    [0, 20, 28, x_{15}],
    [0, 22,  5, 35],
    [0, 23, 37, x_{7}],
    [0, 24, 54, 23],
    [0, 25,  9, x_{1}],
    [0, 28,  7, x_{2}],
    [0, 36, 16, x_{5}],
\\~
    [0, 41, 15, x_{4}],
    [0, 42, 18, x_{6}],
    [0, 43, 25, x_{9}],
    [0, 49, 17, x_{14}]
\\\\
h = 3, n = 19, u = 25 ~(+1~ mod~ 57):
\\~
    [0,  1,  2, x_{25}],
    [0,  2,  4, x_{24}],
    [0,  3,  6, x_{23}],
    [0,  4,  8, x_{22}],
    [0,  5, 10, x_{21}],
    [0,  6, 12, x_{20}],
    [0,  7, 14, x_{19}],
\\~
    [0,  8, 16, x_{18}],
    [0,  9, 24, x_{12}],
    [0, 10, 23, x_{13}],
    [0, 12, 21, x_{14}],
    [0, 13, 25, x_{15}],
    [0, 16, 27, x_{16}],
    [0, 18, 28, x_{17}],
\\~
    [0, 20,  3, 31],
    [0, 23, 37, x_{9}],
    [0, 25,  9, x_{3}],
    [0, 27, 52, x_{2}],
    [0, 31,  1, x_{1}],
    [0, 33, 15, x_{5}],
    [0, 35, 13, x_{6}],
\\~
    [0, 36,  7, x_{4}],
    [0, 40, 17, x_{7}],
    [0, 42, 18, x_{8}],
    [0, 43, 22, x_{10}],
    [0, 46, 26, x_{11}]
\end{array}\]}

\subsubsection*{B3 HSD$(3^{23} u^1)$}

{\footnotesize\[\begin{array}{l}
h = 3, n = 23, u = 17 ~(+1~ mod~ 69):
\\~
    [0,  1,  2, x_{17}],
    [0,  2,  4, x_{16}],
    [0,  3,  6, x_{15}],
    [0,  4,  8, x_{14}],
    [0,  5, 10, x_{13}],
    [0,  7, 24, x_2],
    [0,  8, 20, x_4],
\\~
    [0,  9, 27, 53],
    [0, 10, 64, 31],
    [0, 11, 58, 30],
    [0, 12, 60, 29],
    [0, 13, 62, 27],
    [0, 14, 28, x_3],
    [0, 16, 22, x_{12}],
\\~
    [0, 17, 54, 35],
    [0, 18, 29, x_5],
    [0, 20, 30, x_7],
    [0, 21, 34, x_6],
    [0, 22, 31, x_9],
    [0, 24, 68, 36],
    [0, 25, 33, x_8],
\\~
    [0, 27, 53, 24],
    [0, 30, 37, x_{11}],
    [0, 54, 26, x_1],
    [0, 63, 13, x_{10}]
\\\\
h = 3, n = 23, u = 19 ~(+1~ mod~ 69):
\\~
    [0,  1,  2, x_{19}],
    [0,  2,  4, x_{18}],
    [0,  3,  6, x_{17}],
    [0,  4,  8, x_{16}],
    [0,  5, 10, x_{15}],
    [0,  6, 12, x_{14}],
    [0,  8, 22, x_4],
\\~
    [0,  9, 48, x_2],
    [0, 10, 25, 52],
    [0, 11, 60, 29],
    [0, 12, 56, x_1],
    [0, 13, 64, 28],
    [0, 14, 26, x_5],
    [0, 15, 52, 35],
\\~
    [0, 16, 40, x_6],
    [0, 18, 28, x_7],
    [0, 19, 66, 38],
    [0, 20, 31, x_8],
    [0, 21, 34, x_9],
    [0, 22, 30, x_{11}],
    [0, 24, 33, x_{10}],
\\~
    [0, 25, 68, 36],
    [0, 26, 53, 19],
    [0, 29, 45, x_3],
    [0, 30, 37, x_{13}],
    [0, 62, 14, x_{12}]
\\\\
h = 3, n = 23, u = 23 ~(+1~ mod~ 69):
\\~
    [0,  1,  2, x_{23}],
    [0,  2,  4, x_{22}],
    [0,  3,  6, x_{21}],
    [0,  4,  8, x_{20}],
    [0,  5, 10, x_{19}],
    [0,  6, 12, x_{18}],
    [0,  7, 14, x_{17}],
\\~
    [0,  8, 27, x_5],
    [0, 10, 26, x_8],
    [0, 11, 25, x_9],
    [0, 12, 47, x_4],
    [0, 13, 28, x_7],
    [0, 14, 31, x_6],
    [0, 15, 64, 31],
\\~
    [0, 16, 52, x_3],
    [0, 17, 30, x_{11}],
    [0, 18, 50, x_2],
    [0, 19, 29, x_{12}],
    [0, 20, 62, 40],
    [0, 21, 68, 39],
    [0, 24, 33, x_{14}],
\\~
    [0, 25, 37, x_{13}],
    [0, 26, 34, x_{15}],
    [0, 27, 48, 18],
    [0, 28, 56, 25],
    [0, 32, 58, x_1],
    [0, 34, 45, x_{10}],
    [0, 60, 15, x_{16}]
\end{array}\]}

\subsubsection*{Appendix C1 HSD$(9^4 u^1)$}

{\footnotesize\[\begin{array}{l}
h = 9, n = 4, u = 1 ~(+2~ mod~ 36):
\\~
{\bf [0,  9, 18, 27]},
{\bf [0, 11, 29, 18]},
{\bf [0, 18, 11, 29]},
[0,  2, 17,  7],
[0,  3, 14,  1],
[0,  5, 10, 35],
[0,  6, 25,  3],
\\~
[0,  7, 34, 17],
[0, 10, 21, 23],
[0, 13, 15, x_1],
[0, 14, 27, 33],
[0, 15,  1,  6],
[0, 17,  7, 22],
[0, 27,  6,  5],
\\~[0, 29, 23, 26], [1,  0,  2, x_1]
\\\\
h = 9, n = 4, u = 2 ~(+2~ mod~ 36):
\\~
{\bf [0, 18, 11, 29]},
{\bf [0, 27, 18,  9]},
{\bf [0, 35, 17, 18]},
[0,  1, 26,  3],
[0,  2, 23,  1],
[0,  3, 33, x_1],
[0,  5, 34, 31],
\\~
[0, 11, 22, 17],
[0,  6, 25, x_2],
[0,  7, 21, 34],
[0, 10,  7,  5],
[0, 14,  1, 27],
[0, 15,  5, 22],
[0, 17, 15, 26],
\\~
[0, 21,  6, 35],
[1, 28, 34, x_1], [1, 31, 18, x_2]
\\
\\
h = 9, n = 4, u = 4 ~(+2~ mod~ 36):
\\~
{\bf [0,  1, 18, 19]},
{\bf                  [0, 18,  3, 21]},
{\bf                  [0, 19,  1, 18]},
                  [0,  2,  5, x_4],
                  [0,  3, 10,  1],
                  [0,  5,  7, 26],
                  [0,  6, 11, x_3],
\\~
                  [0,  7, 13, 34],
                  [0,  9, 35, 14],
[0, 10, 17, x_2],
                  [0, 11, 30, 25],
                  [0, 13, 27, 30],
                  [0, 14, 25, x_1],
                  [0, 23, 14, 13],
\\~[0, 25,  2, 31],
                  [1,  3,  4, x_4],
                  [1,  7, 16, x_1],
                  [1, 11, 14, x_3],
                  [1, 23,  8, x_2]
\\
\\
h = 9, n = 4, u = 5 ~(+2~ mod~ 36):
\\~
{\bf    [0,  1, 18, 19]},
{\bf    [0,  9, 27, 18]},
{\bf    [0, 18,  3, 21]},
    [0,  2,  5, x_5],
    [0,  3, 26, 17],
    [0,  5, 19, 26],
    [0,  6, 11, x_4],
\\~
    [0,  7,  6, x_1],
    [0, 10, 17, x_3],
    [0, 11, 22,  9],
    [0, 13, 34, 23],
    [0, 14, 25, x_2],
    [0, 15, 21,  2],
    [0, 17,  7, 22],
\\~
    [0, 31, 29, 30],
    [1,  3,  4, x_5],
    [1,  4, 31, x_1],
    [1, 11, 14, x_4],
    [1, 23, 28, x_3],
    [1, 31,  2, x_2]
\\\\
h = 9, n = 4, u = 7 ~(+2~ mod~ 36):
\\~
{\bf    [0,  1, 19, 18]},
{\bf    [0, 18,  3, 21]},
{\bf    [0, 31, 18, 13]},
    [0,  2,  5, x_7],
    [0,  3, 13, x_3],
    [0,  6, 11, x_6],
    [0,  7,  1, 14],
\\~
    [0,  9, 22, 15],
    [0, 10, 17, x_5],
    [0, 11, 30,  9],
    [0, 13, 27,  6],
    [0, 14, 25, x_4],
    [0, 17, 26, x_1],
    [0, 19, 21, x_2],
\\~
    [0, 27,  2,  1],
    [1,  3,  4, x_7],
    [1,  4,  6, x_2],
    [1, 11, 14, x_6],
    [1, 12, 22, x_3],
    [1, 23, 28, x_5],
    [1, 31,  2, x_4],
\\~
    [1, 32, 15, x_1]
\end{array}\]}

{\footnotesize\[\begin{array}{l}
h = 9, n = 4, u = 8 ~(+2~ mod~ 36):
\\~
{\bf    [0,  1, 18, 19]},
{\bf    [0,  9, 27, 18]},
{\bf    [0, 18,  3, 21]},
    [0,  2,  5, x_8],
    [0,  3, 13, x_4],
    [0,  5, 19, 26],
    [0,  6, 11, x_7],
\\~
    [0,  7,  1, x_1],
    [0, 10, 17, x_6],
    [0, 11, 22, 13],
    [0, 13, 34, 23],
    [0, 14, 25, x_5],
    [0, 17, 26, x_2],
    [0, 19, 21, x_3],
\\~
    [0, 23,  6,  1],
    [1,  2,  8, x_4],
    [1,  3,  4, x_8],
    [1,  4,  6, x_3],
    [1, 11, 14, x_7],
    [1, 16, 30, x_1],
    [1, 22, 31, x_2],
\\~
    [1, 23, 28, x_6],
    [1, 31,  2, x_5]
\\\\
h = 9, n = 4, u = 10 ~(+2~ mod~ 36):
\\~
{\bf    [0,  9, 18, 27]},
{\bf    [0, 17, 35, 18]},
{\bf    [0, 18, 31, 13]},
    [0,  1,  3, x_{10}],
    [0,  2, 11, x_2],
    [0,  3,  2, x_9],
    [0,  6, 25, x_1],
\\~
    [0,  7, 10, x_8],
    [0, 10, 27, 29],
    [0, 13, 14, x_6],
    [0, 14,  1, 31],
    [0, 15,  5, x_5],
    [0, 19, 13, x_4],
    [0, 21,  7, x_3],
\\~
    [0, 31,  6, x_7],
    [1,  2,  4, x_{10}],
    [1,  4,  7, x_8],
    [1,  8, 14, x_5],
    [1, 10, 15, x_7],
    [1, 12, 22, x_4],
    [1, 14, 28, x_3],
\\~
    [1, 23,  8, x_2],
    [1, 26, 11, x_6],
    [1, 27, 16, x_1],
    [1, 32,  3, x_9]
\\\\
h = 9, n = 4, u = 11 ~(+2~ mod~ 36):
\\~
{\bf    [0,  9, 27, 18]},
{\bf    [0, 18,  7, 25]},
{\bf    [0, 31, 18, 13]},
    [0,  1,  3, x_{11}],
    [0,  2,  5, x_{10}],
    [0,  6, 11, x_9],
    [0,  7, 26, x_2],
\\~
    [0, 10, 17, x_8],
    [0, 11, 21, x_6],
    [0, 13, 19, x_5],
    [0, 14, 23, x_7],
    [0, 15, 14, x_1],
    [0, 17,  2, 23],
    [0, 27, 13, x_4],
\\~
    [0, 33,  6, x_3],
    [1,  2,  4, x_{11}],
    [1,  3,  6, x_{10}],
    [1,  8, 27, x_2],
    [1, 11, 12, x_9],
    [1, 12, 23, x_3],
    [1, 14, 28, x_4],
\\~
    [1, 15, 22, x_7],
    [1, 18,  3, x_1],
    [1, 31,  0, x_8],
    [1, 32,  2, x_6],
    [1, 34,  8, x_5]
\\\\
h = 9, n = 4, u = 13 ~(+4~ mod~ 36):
\\~
    [0,  2,  7, x_{13}],
    [0,  7, 10, x_{12}],
    [0,  9, 14, x_{10}],
    [0, 11, 18, x_8],
    [0, 14, 17, x_1],
    [0, 17, 27, x_5],
    [0, 19, 33, x_4],
\\~
    [0, 21, 30, x_6],
    [0, 22,  5, x_2],
    [0, 23, 22,  1],
    [0, 25, 15, x_7],
    [0, 27, 25, x_3],
    [0, 29,  6, x_{11}],
    [0, 30, 11, x_9],
\\~
    [1,  0,  7, x_{11}],
    [1,  3,  0, x_1],
    [1,  4, 10, x_{13}],
    [1,  6,  8, x_{12}],
    [1,  7, 18, x_7],
    [1, 10, 19, x_6],
    [1, 11, 24, x_2],
\\~
    [1, 18, 28, x_4],
    [1, 23,  2, x_5],
    [1, 24, 15, x_8],
    [1, 30, 16, x_3],
    [1, 32,  6, x_9],
    [1, 34, 11, x_{10}],
    [2,  1,  7, x_{12}],
\\~
    [2,  4, 25, x_5],
    [2,  7,  9, x_{13}],
    [2, 11, 17, x_{11}],
    [2, 12, 27, x_4],
    [2, 13, 31, x_3],
    [2, 15,  4, x_{10}],
    [2, 19,  5, x_9],
\\~
    [2, 20, 19, x_1],
    [2, 25,  0, x_8],
    [2, 28, 23, x_2],
    [2, 31, 12, x_6],
    [2, 32, 13, x_7],
    [3,  0,  5, x_{12}],
    [3,  2,  9, x_8],
\\~
    [3,  4, 13, x_{10}],
    [3,  5,  8, x_{13}],
    [3,  6,  4, x_{11}],
    [3,  8, 30, x_3],
    [3,  9, 16, x_9],
    [3, 13,  2, x_2],
    [3, 14, 20, x_7],
\\~
    [3, 18,  0, x_5],
    [3, 21,  6, x_4],
    [3, 24,  1, x_6],
    [3, 25, 26, x_1]
\end{array}\]}

\subsubsection*{C2 HSD$(9^5 u^1)$}

{\footnotesize\[\begin{array}{l}
h = 9, n = 5, u = 2 ~(+1~ mod~ 45):
\\~
    [0,  1,  2, x_{8}],
    [0,  2, 34, 18],
    [0,  4, 41, 23],
    [0,  6,  8, x_{6}],
    [0,  7, 36, 19],
    [0,  9, 18, x_{3}],
    [0, 11, 14, x_{5}],
\\~
    [0, 12, 24, 11],
    [0, 14,  7, 31],
    [0, 19, 23, x_{4}],
    [0, 22, 39, x_{1}],
    [0, 37, 13, x_{2}],
    [0, 42,  3, x_{7}]
\\\\
n = 5, u = 4 ~(+1~ mod~ 45):
\\~[0,  1,  2, x_4],
                  [0,  2, 29,  8],
                  [0,  4, 21, 38],
                  [0,  6,  8, x_2],
                  [0,  7, 36, 24],
                  [0,  8, 19, 41],
                  [0,  9, 22, 36],
\\~[0, 11, 14, x_1],
                  [0, 13, 32, 14],
                  [0, 16,  4, 23],
                  [0, 42,  3, x_3]
\\\\
h = 9, n = 5, u = 8 ~(+1~ mod~ 45):
\\~
    [0,  1,  2, x_{8}],
    [0,  2, 34, 18],
    [0,  4, 41, 23],
    [0,  6,  8, x_{6}],
    [0,  7, 36, 19],
    [0,  9, 18, x_{3}],
    [0, 11, 14, x_{5}],
\\~
    [0, 12, 24, 11],
    [0, 14,  7, 31],
    [0, 19, 23, x_{4}],
    [0, 22, 39, x_{1}],
    [0, 37, 13, x_{2}],
    [0, 42,  3, x_{7}]
\\\\
h = 9, n = 5, u = 10 ~(+1~ mod~ 45):
\\~[0,  1, 17,  4], [0,  3, 11, x_6], [0,  6, 38, x_7], [0,  7, 14, x_3],
[0,  8, 27, x_2], [0, 11, 39, 27], [0, 36, 12, x_5],
\\~[0,  2, 24,  6], [0, 14, 23,  1], [0, 17, 19, x_{10}], [0, 21,  9, x_4],
[0, 26, 37, x_1], [0, 29, 43, x_9], [0, 41, 44, x_8]
\\\\
h = 9, n = 5, u = 14 ~(+1~ mod~ 45):
\\~
    [0,  1,  2, x_{14}],
    [0,  2,  4, x_{13}],
    [0,  3,  6, x_{12}],
    [0,  6, 14, x_{8}],
    [0,  7, 34, x_{5}],
    [0,  9, 18, x_{4}],
    [0, 11,  7, 33],
\\~
    [0, 12, 28, x_{6}],
    [0, 13, 19, x_{10}],
    [0, 14, 42, 26],
    [0, 17, 24, x_{9}],
    [0, 18, 29, x_{7}],
    [0, 21, 44, x_{2}],
    [0, 22,  8, x_{1}],
\\~
    [0, 37, 13, x_{3}],
    [0, 41,  9, x_{11}]
\end{array}\]}

{\footnotesize\[\begin{array}{l}
h = 9, n = 5, u = 16 ~(+1~ mod~ 45):
\\~
    [0,  1,  2, x_{16}],
    [0,  2,  4, x_{15}],
    [0,  3,  6, x_{14}],
    [0,  4,  8, x_{13}],
    [0,  6, 14, x_{10}],
    [0,  7, 28, x_{6}],
    [0,  8, 22, x_{7}],
\\~
    [0,  9, 18, x_{8}],
    [0, 11, 33, x_{4}],
    [0, 12, 38, x_{5}],
    [0, 13, 19, x_{12}],
    [0, 14, 32, x_{3}],
    [0, 16, 44, x_{2}],
    [0, 17, 24, x_{11}],
\\~
    [0, 18, 29, x_{9}],
    [0, 21,  9, 32],
    [0, 26, 42, x_{1}]
\end{array}\]}

\subsubsection*{C3 HSD$(9^9 u^1)$}

{\footnotesize\[\begin{array}{l}
h = n = 9, u = 8 ~(+1~ mod~ 81):
\\~
    [0,  1,  2, x_{8}],
    [0,  2, 69, 49],
    [0,  3, 11, x_{2}],
    [0,  4,  6, x_{7}],
    [0,  5, 74, 57],
    [0,  6, 31, 23],
    [0,  7, 71, 55],
\\~
    [0, 10, 20, x_{1}],
    [0, 11, 39, 51],
    [0, 13, 16, x_{6}],
    [0, 14, 58, 11],
    [0, 15, 53, 29],
    [0, 19, 66, 40],
    [0, 21,  5, 46],
\\~
    [0, 22, 26, x_{5}],
    [0, 23, 73, 22],
    [0, 25, 46, 74],
    [0, 29, 68, 33],
    [0, 31, 37, x_{4}],
    [0, 32, 51, 13],
    [0, 33, 38, x_{3}],
\\~
    [0, 37, 57, 15]
\\\\
h = n = 9, u = 10 ~(+1~ mod~ 81):
\\~
    [0,  1,  2, x_{10}],
    [0,  2, 66, 50],
    [0,  3, 11, x_{4}],
    [0,  4,  6, x_{9}],
    [0,  5, 67, 74],
    [0,  6, 74, 48],
    [0,  8, 22, 57],
\\~
    [0, 10, 20, x_{3}],
    [0, 11, 53, 32],
    [0, 12, 71, 47],
    [0, 13, 16, x_{8}],
    [0, 14, 40, 65],
    [0, 15, 68, 38],
    [0, 17, 29, x_{2}],
\\~
    [0, 19,  4, 44],
    [0, 20, 80, 51],
    [0, 22, 26, x_{7}],
    [0, 23, 34, x_{1}],
    [0, 28,  5, 52],
    [0, 31, 37, x_{6}],
\\~
    [0, 32, 78, 40],
    [0, 33, 38, x_{5}],
    [0, 37, 62, 20]
\\\\
h = n = 9, u = 14 ~(+1~ mod~ 81):
\\~
    [0,  1,  2, x_{14}],
    [0,  2,  4, x_{13}],
    [0,  3, 24, 31],
    [0,  4, 68, 62],
    [0,  5, 30, 61],
    [0,  8, 12, x_{10}],
    [0, 10, 21, x_{3}],
\\~
    [0, 11, 78, 55],
    [0, 12, 15, x_{12}],
    [0, 13, 64, 48],
    [0, 14, 61, 40],
    [0, 15, 31, x_{2}],
    [0, 17, 22, x_{11}],
    [0, 19, 29, x_{5}],
\\~
    [0, 20, 32, x_{1}],
    [0, 22, 75, 38],
    [0, 24, 39, x_{4}],
    [0, 25, 38, x_{6}],
    [0, 26, 33, x_{8}],
    [0, 28, 34, x_{9}],
    [0, 29, 73, 22],
\\~
    [0, 32, 40, x_{7}],
    [0, 33, 62, 23],
    [0, 34, 80, 39],
    [0, 35, 67, 24]
\\\\
h = n = 9, u = 16 ~(+1~ mod~ 81):
\\~
    [0,  1,  2, x_{16}],
    [0,  2,  4, x_{15}],
    [0,  3,  6, x_{14}],
    [0,  4, 53, 48],
    [0,  6, 65, 30],
    [0,  7, 20, x_{1}],
    [0,  8, 13, x_{13}],
\\~
    [0, 10, 22, x_{5}],
    [0, 11, 71, 57],
    [0, 12, 64, 31],
    [0, 13, 30, x_{2}],
    [0, 15, 70, 41],
    [0, 16, 31, x_{6}],
    [0, 17, 21, x_{12}],
\\~
    [0, 19, 29, x_{8}],
    [0, 20, 73, 43],
    [0, 21, 32, x_{7}],
    [0, 22, 38, x_{3}],
    [0, 23,  3, 37],
    [0, 24, 66, 23],
    [0, 25, 39, x_{4}],
\\~
    [0, 26, 33, x_{10}],
    [0, 28, 34, x_{11}],
    [0, 31, 12, 56],
    [0, 32, 40, x_{9}],
    [0, 39, 74, 34]
\end{array}\]}

\subsubsection*{Appendix D HSD$(4^n u^1)$ for $n=19, 22$}

{\footnotesize\[\begin{array}{l}
h = 4, n = 19, u = 31 ~(+2~ mod~ 76):
\\~
    [0,  2,  4, x_{10}],
    [0,  3,  7, x_{18}],
    [0,  4,  9, x_{30}],
    [0,  5,  6, x_{31}],
    [0,  6, 13, x_{28}],
    [0,  9, 17, x_{27}],
    [0, 10, 24, x_{20}],
\\~
    [0, 11, 20, x_{24}],
    [0, 12, 21, x_{26}],
    [0, 13, 12, x_{5}],
    [0, 14, 36, x_{16}],
    [0, 15, 27, x_{22}],
    [0, 16, 34, x_{17}],
    [0, 17, 30, x_{21}],
\\~
    [0, 20, 55, x_{11}],
    [0, 21, 37, x_{19}],
    [0, 23, 44, x_{15}],
    [0, 25, 47, x_{14}],
    [0, 27, 54, x_{12}],
    [0, 28, 60, x_{7}],
    [0, 30, 58, 45],
\\~
    [0, 31, 68, 34],
    [0, 32, 61, x_{13}],
    [0, 35, 67, x_{6}],
    [0, 36, 66, 37],
    [0, 39, 50, x_{23}],
    [0, 41, 11, x_{1}],
    [0, 50, 14, x_{9}],
\\~
    [0, 52, 23, x_{4}],
    [0, 53,  5, 41],
    [0, 54, 29, x_{3}],
    [0, 58, 25, x_{8}],
    [0, 61, 28, x_{2}],
    [0, 68,  2, x_{25}],
    [0, 73,  3, x_{29}],
\\~
    [1,  0,  3, x_{31}],
    [1,  2,  6, x_{18}],
    [1,  3,  5, x_{10}],
    [1,  5,  8, x_{30}],
    [1,  6, 12, x_{29}],
    [1,  7, 17, x_{25}],
    [1,  8, 19, x_{24}],
\\~
    [1,  9, 14, x_{28}],
    [1, 10, 71, x_{23}],
    [1, 11, 18, x_{26}],
    [1, 12, 32, x_{14}],
    [1, 13, 31, x_{17}],
    [1, 15, 35, x_{16}],
    [1, 18, 44, x_{6}],
\\~
    [1, 19, 42, x_{8}],
    [1, 22, 38, x_{19}],
    [1, 25, 49, x_{9}],
    [1, 26, 41, x_{15}],
    [1, 27, 69, 37],
    [1, 28, 45, x_{12}],
    [1, 29, 55, x_{7}],
\\~
    [1, 32, 63, x_{2}],
    [1, 34, 46, x_{22}],
    [1, 40, 53, x_{21}],
    [1, 43, 57, x_{20}],
    [1, 44, 27, x_{5}],
    [1, 47, 26, x_{3}],
    [1, 48, 24, x_{1}],
\\~
    [1, 55, 28, x_{4}],
    [1, 57, 34, x_{11}],
    [1, 61, 36, x_{13}],
    [1, 70,  2, x_{27}],
\end{array}\]}

{\footnotesize\[\begin{array}{l}
h = 4, n = 19, u = 30 ~(+1~ mod~ 76):
\\~
    [0,  1,  2, x_{30}],
    [0,  2,  4, x_{29}],
    [0,  3,  6, x_{28}],
    [0,  4,  8, x_{27}],
    [0,  5, 10, x_{26}],
    [0,  6, 12, x_{25}],
    [0,  7, 14, x_{24}],
\\~
    [0,  8, 16, x_{23}],
    [0,  9, 18, x_{22}],
    [0, 10, 20, x_{21}],
    [0, 12, 25, x_{17}],
    [0, 16, 33, x_{11}],
    [0, 17, 35, x_{10}],
    [0, 20, 36, x_{12}],
\\~
    [0, 21,  1, x_{1}],
    [0, 22, 37, x_{14}],
    [0, 26, 71, 41],
    [0, 27,  3, x_{3}],
    [0, 28,  7, 39],
    [0, 29, 52, x_{2}],
    [0, 31, 42, x_{20}],
\\~
    [0, 34, 46, x_{18}],
    [0, 35, 63, 26],
    [0, 36, 50, x_{16}],
    [0, 43, 21, x_{4}],
    [0, 51, 17, x_{5}],
    [0, 52, 23, x_{6}],
    [0, 53, 28, x_{7}],
\\~
    [0, 58, 31, x_{8}],
    [0, 61, 29, x_{13}],
    [0, 62, 32, x_{9}],
    [0, 63, 27, x_{15}],
    [0, 65, 22, x_{19}],
\\\\
h = 4, n = 19, u = 33 ~(+2~ mod~ 76):
\\~
    [0,  1,  3, x_{10}],
    [0,  2,  5, x_{1}],
    [0,  4,  6, x_{33}],
    [0,  5,  8, x_{32}],
    [0,  6, 12, x_{31}],
    [0,  7, 14, x_{28}],
    [0,  8, 15, x_{30}],
\\~
    [0,  9, 17, x_{29}],
    [0, 10, 20, x_{27}],
    [0, 11, 25, x_{19}],
    [0, 12, 23, x_{26}],
    [0, 14, 28, x_{24}],
    [0, 15, 37, x_{17}],
    [0, 16, 29, x_{22}],
\\~
    [0, 18, 52, x_{12}],
    [0, 20, 43, x_{18}],
    [0, 21, 34, x_{21}],
    [0, 22, 51, x_{6}],
    [0, 24, 44, x_{13}],
    [0, 25, 40, x_{20}],
    [0, 26, 59, x_{3}],
\\~
    [0, 27, 47, x_{15}],
    [0, 28, 58, x_{7}],
    [0, 29, 50, x_{11}],
    [0, 30, 61, x_{2}],
    [0, 31, 54, x_{8}],
    [0, 32, 72, x_{9}],
    [0, 33, 60, x_{23}],
\\~
    [0, 35, 74, 45],
    [0, 36,  1, 47],
    [0, 37, 49, x_{25}],
    [0, 41,  7, 32],
    [0, 42, 66, x_{4}],
    [0, 49, 13, x_{5}],
    [0, 61, 35, x_{14}],
\\~
    [0, 69, 30, x_{16}],
    [1,  2,  7, x_{32}],
    [1,  3,  4, x_{1}],
    [1,  4, 12, x_{29}],
    [1,  5,  9, x_{33}],
    [1,  6, 15, x_{28}],
    [1,  7, 13, x_{31}],
\\~
    [1,  9, 14, x_{30}],
    [1, 10, 35, x_{16}],
    [1, 11, 21, x_{27}],
    [1, 12, 24, x_{25}],
    [1, 13, 22, x_{26}],
    [1, 14, 32, x_{15}],
    [1, 15, 31, x_{13}],
\\~
    [1, 18, 44, x_{14}],
    [1, 19, 36, x_{18}],
    [1, 22, 38, x_{19}],
    [1, 24, 51, x_{8}],
    [1, 25, 50, x_{6}],
    [1, 29, 61, x_{4}],
    [1, 32, 33, x_{11}],
\\~
    [1, 34, 56, x_{5}],
    [1, 35, 46, x_{22}],
    [1, 37, 19, x_{24}],
    [1, 38, 55, x_{23}],
    [1, 45, 10, x_{2}],
    [1, 51,  8, x_{3}],
    [1, 54,  6, x_{17}],
\\~
    [1, 55, 25, x_{7}],
    [1, 57, 29, x_{9}],
    [1, 60,  5, x_{20}],
    [1, 61, 37, x_{12}],
    [1, 64,  3, x_{21}],
    [1, 74,  2, x_{10}],
\\\\
h = 4, n = 19, u = 34 ~(+1~ mod~ 76):
\\~
    [0,  1,  2, x_{34}],
    [0,  2,  4, x_{33}],
    [0,  3,  6, x_{32}],
    [0,  4,  8, x_{31}],
    [0,  5, 10, x_{30}],
    [0,  6, 12, x_{29}],
    [0,  7, 14, x_{28}],
\\~
    [0,  8, 16, x_{27}],
    [0,  9, 18, x_{26}],
    [0, 10, 20, x_{25}],
    [0, 11, 22, x_{24}],
    [0, 12, 26, x_{19}],
    [0, 14, 35, x_{11}],
    [0, 17, 34, x_{12}],
\\~
    [0, 20, 36, x_{15}],
    [0, 22, 37, x_{18}],
    [0, 24, 51, x_{8}],
    [0, 25,  7, x_{2}],
    [0, 26, 63, x_{7}],
    [0, 28, 53, x_{6}],
    [0, 29,  3, 34],
\\~
    [0, 30, 43, x_{20}],
    [0, 32, 44, x_{21}],
    [0, 35, 15, x_{4}],
    [0, 37,  1, x_{3}],
    [0, 40, 17, x_{22}],
    [0, 42,  9, x_{23}],
    [0, 43, 65, x_{1}],
\\~
    [0, 49, 21, x_{5}],
    [0, 53, 24, x_{9}],
    [0, 55, 31, x_{10}],
    [0, 58, 27, x_{14}],
    [0, 60, 30, x_{13}],
    [0, 61, 29, x_{16}],
    [0, 63, 28, x_{17}],
\\\\
h = 4, n = 19, u = 35 ~(+2~ mod~ 76):
\\~
    [0,  2,  6, x_{23}],
    [0,  3,  5, x_{26}],
    [0,  4,  7, x_{35}],
    [0,  5,  8, x_{34}],
    [0,  6, 13, x_{32}],
    [0,  8, 24, x_{13}],
    [0,  9, 17, x_{31}],
\\~
    [0, 11, 20, x_{25}],
    [0, 12, 21, x_{30}],
    [0, 14, 34, x_{12}],
    [0, 15, 39, x_{21}],
    [0, 17, 28, x_{17}],
    [0, 18, 47, x_{14}],
    [0, 22, 35, x_{20}],
\\~
    [0, 23, 37, x_{19}],
    [0, 24, 41, x_{16}],
    [0, 27, 45, x_{15}],
    [0, 28, 43, x_{22}],
    [0, 29,  2, x_{27}],
    [0, 30, 12, x_{28}],
    [0, 31, 62, x_{1}],
\\~
    [0, 33, 66, 35],
    [0, 34, 46, x_{24}],
    [0, 35, 50, x_{8}],
    [0, 40, 11, x_{3}],
    [0, 41,  1, x_{2}],
    [0, 43, 60, x_{5}],
    [0, 44, 54, x_{29}],
\\~
    [0, 50, 18, x_{4}],
    [0, 51, 29, x_{7}],
    [0, 53, 25, x_{6}],
    [0, 56, 33, x_{9}],
    [0, 60, 36, x_{11}],
    [0, 65, 32, x_{10}],
    [0, 66,  4, x_{18}],
\\~
    [0, 73,  3, x_{33}],
    [1,  0,  5, x_{34}],
    [1,  2,  4, x_{26}],
    [1,  3,  7, x_{23}],
    [1,  5,  6, x_{35}],
    [1,  6, 12, x_{33}],
    [1,  7, 17, x_{29}],
\\~
    [1,  8, 45, x_{27}],
    [1,  9, 14, x_{32}],
    [1, 10, 21, x_{25}],
    [1, 11, 18, x_{30}],
    [1, 13, 25, x_{24}],
    [1, 14, 55, x_{5}],
    [1, 16, 37, x_{17}],
\\~
    [1, 17, 16, x_{20}],
    [1, 18, 54, x_{6}],
    [1, 19, 35, x_{18}],
    [1, 21, 51, x_{12}],
    [1, 22, 50, x_{7}],
    [1, 23, 49, x_{28}],
    [1, 28, 59, x_{10}],
\\~
    [1, 29, 63, x_{11}],
    [1, 30, 69, x_{8}],
    [1, 33, 46, x_{9}],
    [1, 35, 67, x_{4}],
    [1, 37, 62, x_{3}],
    [1, 38, 68, x_{2}],
    [1, 40, 13, x_{1}],
\\~
    [1, 47, 24, x_{16}],
    [1, 51, 31, x_{13}],
    [1, 52, 26, x_{15}],
    [1, 53, 28, x_{14}],
    [1, 56, 22, x_{19}],
    [1, 63,  8, x_{22}],
    [1, 64, 10, x_{21}],
\\~
    [1, 70,  2, x_{31}],
\\\\
h = 4, n = 22, u = 34 ~(+1~ mod~ 88):
\\~
    [0,  1,  2, x_{34}],
    [0,  2,  4, x_{33}],
    [0,  3,  6, x_{32}],
    [0,  4,  8, x_{31}],
    [0,  5, 10, x_{30}],
    [0,  6, 12, x_{29}],
    [0,  7, 14, x_{28}],
\\~
    [0,  8, 16, x_{27}],
    [0,  9, 18, x_{26}],
    [0, 10, 20, x_{25}],
    [0, 11, 26, x_{21}],
    [0, 12, 23, x_{24}],
    [0, 14, 28, x_{20}],
    [0, 16, 37, x_{14}],
\\~
    [0, 17, 40, x_{15}],
    [0, 19, 35, x_{19}],
    [0, 20, 39, x_{16}],
    [0, 24, 42, x_{17}],
    [0, 25, 67, 20],
    [0, 26, 43, x_{18}],
    [0, 32, 77, 31],
\\~
    [0, 33,  9, x_{3}],
    [0, 34, 63, 27],
    [0, 35, 47, x_{23}],
    [0, 37, 50, x_{22}],
    [0, 39, 71, 26],
    [0, 40, 15, x_{1}],
    [0, 50,  3, x_{2}],
\\~
    [0, 57, 19, x_{4}],
    [0, 58, 30, x_{5}],
    [0, 59, 29, x_{6}],
    [0, 60, 27, x_{9}],
    [0, 61, 24, x_{7}],
    [0, 65, 31, x_{10}],
    [0, 67, 32, x_{8}],
\\~
    [0, 70, 34, x_{12}],
    [0, 73, 33, x_{11}],
    [0, 75, 36, x_{13}]
\end{array}\]}

\end{document}